\documentclass[]{scrartcl}

\setkomafont{disposition}{\normalfont\bfseries}
\usepackage{amsmath}
\usepackage{amsfonts}
\usepackage{graphicx}
\usepackage[none]{hyphenat}
\usepackage{float}
\usepackage{fullpage}
\usepackage{calc}
\newcommand{\explain}[2]{\underbrace{#1}_\text{#2}}
\usepackage{multirow}
\usepackage{caption}
\usepackage{subcaption}
\usepackage{wrapfig}

\author{Manan Gandhi}
\title{Trajectory Optimization Algorithm Studies}
\subtitle{A Comparison between DDP and Pseudospectral Methods}

\begin{document}
\maketitle

\begin{abstract}
In complex engineered systems, completing an objective is sometimes not enough. The system must be able to reach a set performance characteristic, such as an unmanned aerial vehicle flying from point A to point B, \textit{under 10 seconds}. This introduces the notion of optimality, what is the most efficient, the fastest, the cheapest way to complete a task. This report explores the two pillars of optimal control, Bellman's Dynamic Programming and Pontryagin's Maximum Principle, and compares implementations of both theories onto simulated systems. Dynamic Programming is realized through a Differential Dynamic Programming Algorithm, where utilizes a forward-backward pass to iteratively optimize a control sequence and trajectory. The Maximum Principle is realized via Gauss Pseudospectral Optimal Control, where the optimal control problem is first approximated through polynomial basis functions, then solved, with optimality being achieved through the costate equations of the Maximum Principle. The results of the report show that, for short time Horizons, DDP can optimize quickly and can generate a trajectory that utilizes less control effort for the same problem formulation. On the other hand Pseudospectral methods can optimize faster for longer time horizons, but require a key understanding of the problem structure. Future work involves completing an implementation of DDP in a C++ code, and testing the speed of convergence for both methods, as well as extended the Pseudospectral Optimal Control framework in to the world of stochastic optimal control.
\end{abstract}

\section{Introduction}
As systems become more complex and tasks become more demanding, optimal control of multi-body systems becomes increasingly challenging. There are many methods built for various classes of optimal control problems, but this paper aims to compare two specific methods, namely Differential Dynamic Programming, and Gauss Pseudospectral Optimal Control in order to compare their respective ability to handle complex dynamic and algebraic constraints. The impact of this project can enter the realm of control with contact dynamics, a notoriously difficult problem with automation, as well as stochastic control. In order to discuss the aforementioned methods for optimal control, we first understand their roots in the two pillars of optimal control, Bellman's Dynamic Programming, and Pontryagin's Maximum Principle.

\subsection{Dynamic Programming}
The basis of Dynamic Programming lies in discretizing the path from the initial state to the final state, then finding the optimal path for each discrete step. This approach will then compare the most optimal path through the entirety of the states, and compute a global optimum. Dynamic Programming provides a sufficient condition for an optimal trajectory. The primary weakness, however, is that it must compute the optimal trajectory for every possible discrete interval. As a result, it suffers from the curse of dimensionality as the number of states and the resolution of each step increases. A consequence of the theory of dynamic programming leads to the Hamilton-Jacobi-Bellman Equation.

\subsubsection*{Hamilton-Jacobi-Bellman Equation}
\begin{eqnarray}
\mathbf{\dot{x}} (t)= \mathbf{f}[\mathbf{x}(t),\mathbf{u}(t),t]
\\V[\mathbf{x}(t),t]=\underset{\mathbf{u}}{\text{min}} \bigg[ \int_{t}^{t+\Delta t}L[\mathbf{x}(t),\mathbf{u}(t),t]dt+V[\mathbf{x'}(t),t+\Delta t]\bigg]
\end{eqnarray}

\begin{itemize}
\item $V[\mathbf{x}(t),t]$ is the cost to go to state at time $t$ in the system. 
\item $V[\mathbf{x'}(t),t+\Delta t]$ is the penalty or the cost to go to the next state after time $\Delta t$.
\item $L[\mathbf{x}(t),\mathbf{u}(t),t]$ is a running cost that accrues as control is utilized to move from one state to another.
\item $\mathbf{x'}(t)=\mathbf{x}(t+\Delta t)$
\end{itemize}

The next step is to be able to solve this cost function $V$ from one state to another over time, and to minimize the control $\mathbf{u}(t)$ over that time interval.

We can add and subtract $V[[\mathbf{x}(t),t]]$ on the right-hand-side of the equation.
\begin{eqnarray}
V[\mathbf{x}(t),t]=\underset{\mathbf{u}}{\text{min}} \bigg[ \int_{t}^{t+\Delta t}L[\mathbf{x}(t),\mathbf{u}(t),t]dt+V[\mathbf{x'}(t),t+\Delta t]-V[\mathbf{x}(t),t]+V[\mathbf{x}(t),t]\bigg]
\\ \nonumber
\\dV= V[\mathbf{x'}(t),t+\Delta t]-V[\mathbf{x}(t),t]
\\ \nonumber
\\V[\mathbf{x}(t),t]=\underset{\mathbf{u}}{\text{min}} \bigg[ \int_{t}^{t+\Delta t}L[\mathbf{x}(t),\mathbf{u}(t),t]dt+dV+V[\mathbf{x}(t),t]\bigg] \label{5}
\end{eqnarray}

Next we can utilize Taylor Polynomials to create an expression for $dV$ that can be substituted into Equation \eqref{5}. At this point the notation expressing $\mathbf{x}$ or $\mathbf{u}$ as a function of time will be dropped: $\mathbf{x}(t) = \mathbf{x}$.

\begin{eqnarray}
V[x',t+\Delta t] =&& V[x'-x+x,t+\Delta t] = V[x+\delta x,t+\Delta t]
\\=&& V[x,t]+\frac{\partial V}{\partial t}\delta t + V_x^T \delta x 
+ \frac{1}{2} \begin{bmatrix} \delta x\\ \delta t \end{bmatrix}^T
\begin{bmatrix} V_x{_x} & V_x{_t} \\ V_t{_x} & V_t{_t} \end{bmatrix} \begin{bmatrix}\delta x\\ \delta t\end{bmatrix} \label{7}
\end{eqnarray}

With this formulation of $dV$ utilize the fact $\delta \mathbf{x} = f[\mathbf{x},\mathbf{u},t] \delta t$ and the fact that higher order terms go to zero as they are very small. Substituting this in Equation \eqref{7} gives

\begin{eqnarray}
V[x',t+\Delta t] = V[x,t]+\frac{\partial V}{\partial t}\delta t + V_x^T f[\mathbf{x},\mathbf{u},t] \delta t
+ \frac{1}{2} \begin{bmatrix} f[\mathbf{x},\mathbf{u},t] \delta t\\ \delta t \end{bmatrix}^T
\begin{bmatrix} V_x{_x} & V_x{_t} \\ V_t{_x} & V_t{_t} \end{bmatrix} \begin{bmatrix}f[\mathbf{x},\mathbf{u},t] \delta t\\ \delta t\end{bmatrix}
\\\nonumber
\\\frac{1}{2} \begin{bmatrix} f[\mathbf{x},\mathbf{u},t] \delta t\\ \delta t \end{bmatrix}^T
\begin{bmatrix} V_x{_x} & V_x{_t} \\ V_t{_x} & V_t{_t} \end{bmatrix} \begin{bmatrix}f[\mathbf{x},\mathbf{u},t] \delta t\\ \delta t\end{bmatrix} = 0  \label{9}
\end{eqnarray}

Equation \eqref{9} is true because higher order terms go to zero. The factor of $\delta t$ in each term of the matrices ensures that there will be a $\delta t^2$ in every piece of Equation \eqref{9}.
Next, bringing over $V[x,t]$ over to the left-hand-side will give:

\begin{eqnarray}
dV = \frac{\partial V}{\partial t}\delta t + V_x^T f[\mathbf{x},\mathbf{u},t] \delta t
\\
V[\mathbf{x},t]=\underset{\mathbf{u}}{\text{min}} \bigg[ \int_{t}^{t+\Delta t}L[\mathbf{x},\mathbf{u}),t]dt+\frac{\partial V}{\partial t}\delta t + V_x^T f[\mathbf{x},\mathbf{u},t] \delta t+V[\mathbf{x},t]\bigg]
\end{eqnarray}

Cancelling out $V[\mathbf{x}(t),t]$ on both sides and applying $\int_{t}^{t+\Delta t}L[\mathbf{x},\mathbf{u},t]dt = L[\mathbf{x},\mathbf{u},t] \delta t$ allows $\delta t$ to be divided out of the equation. Finally move $\frac{\partial V}{\partial t}$ to the left hand side and the result is Equation \eqref{final}.

\begin{eqnarray}
0 =&& \underset{\mathbf{u}}{\text{min}} \bigg[ \int_{t}^{t+\Delta t}L[\mathbf{x},\mathbf{u},t]dt+\frac{\partial V}{\partial t}\delta t + V_x^T f[\mathbf{x},\mathbf{u},t] \delta t\bigg]
\\=&&\underset{\mathbf{u}}{\text{min}} \bigg[ L[\mathbf{x},\mathbf{u},t] \delta t+\frac{\partial V}{\partial t}\delta t + V_x^T f[\mathbf{x},\mathbf{u},t] \delta t\bigg]
\\=&&\underset{\mathbf{u}}{\text{min}} \bigg[ L[\mathbf{x},\mathbf{u},t] +\frac{\partial V}{\partial t} + V_x^T f[\mathbf{x},\mathbf{u},t]\bigg]
\\\frac{\partial V}{\partial t} =&& \underset{\mathbf{u}}{\text{min}} \bigg[ L[\mathbf{x},\mathbf{u},t]+V_x^T f[\mathbf{x},\mathbf{u},t]\bigg] \label{final}
\end{eqnarray} 

Thus Equation \eqref{final} gives the Hamilton-Jacobi-Bellman Equation. This equation provides a basis to solve for a cost to go from one end of the state to another based on the running cost $L$ and the terminal cost at the final state.

Now that we have the equation, it is necessary to determine the boundary condition for this non-linear partial differential equation. If we are trying to move from state B to state C, the total cost of going to state C will be the terminal cost at C, plus any sort of running cost that builds up as we input control into the system. This is shown in Equation \eqref{16}.

\begin{eqnarray}
V[\mathbf{x}(t_B),t_B]=\underset{\mathbf{u}}{\text{min}} \bigg[ \int_{t_B}^{t_C} L[\mathbf{x}(t),\mathbf{u}(t),t]dt+V[\mathbf{x'}(t_C),t_C]\bigg] \label{16}
\end{eqnarray}

As a result, in the HJB equation, we can solve for the value function utilizing the "cost to go" or the terminal cost of state C as a boundary condition. This is where the terminal boundary condition is utilized in dynamic programming. Once the value function has been solved over time, we can find the "cost to go" required to be at state B. Now that this is known, we can repeat the process for a prior state movement from A to B.

\subsection{The Maximum Principle}
Pontryagin's Maximum Principle can be stated as necessary conditions for optimality through the use of calculus of variations. An optimal control history will cause the cost function to remain stationary with respect to the control function within the constraint of the system dynamics and initial state ~\cite{Stengel}.

This approach is general and does not require the specification of the shape of the control function before calculation, unlike parametric optimization. The controls can even include piecewise continuous functions.

A general non-linear, time-varying dynamic system can be defined with the 
following equation. And then the derivative of the state can be brought to the right hand side of the equation creating an equality constraint that is represented by the state dynamics. When this equation is satisfied, the equality constraint is zero.

\begin{eqnarray}
\mathbf{\dot{x}} (t)= \mathbf{f}[\mathbf{x}(t),\mathbf{u}(t),t]
\\\mathbf{f}[\mathbf{x}(t),\mathbf{u}(t),t]-\mathbf{\dot{x}} (t)= 0 \label{2}
\end{eqnarray}

Now, the constraint can be adjoined to the cost function in order to create a single equation that can be minimized with respect to the control. Equation \eqref{2} can be combined into the cost function J via a vector of Lagrange multipliers.

\begin{eqnarray}
J= \mathbf{\phi}[\mathbf{x}(t_f),t_f]+\int_{t_0}^{t_f}\mathcal{L}[[\mathbf{x}(t),\mathbf{u}(t),t]+\lambda^T[\mathbf{f}[\mathbf{x}(t),\mathbf{u}(t),t]-\mathbf{\dot{x}} (t)]]dt
\end{eqnarray}

\begin{itemize}
  \item $\mathcal{L}$ is the Lagrangian of the system, the term represents the energy required to change the state from the initial to the final. 
  \item $\mathbf{\phi}$ is the penalty or the cost at the final state.
  \item $\lambda$ is the vector of Lagrange multipliers that adjoins the dynamic constraint to the cost function.

\end{itemize}

The next step is to utilize the Hamiltonian, $\mathcal{H}$ to simplify the equations to manipulate them further via integration by parts. 

\begin{eqnarray}
\mathcal{H}[\mathbf{x}(t),\mathbf{u}(t),\lambda(t),t]= \mathcal{L}[\mathbf{x}(t),\mathbf{u}(t),t]+\lambda^T\mathbf{f}[\mathbf{x}(t),\mathbf{u}(t),t]
\\J= \mathbf{\phi}[\mathbf{x}(t_f),t_f]+\int_{t_0}^{t_f}[\mathcal{H}[\mathbf{x}(t),\mathbf{u}(t),\lambda(t),t]-\lambda^T\mathbf{\dot{x}}(t)]dt
\end{eqnarray}
\begin{eqnarray}
J= \\&&\mathbf{\phi}[\mathbf{x}(t_f),t_f]+[\lambda^T(t_o)\mathbf{x}(t_0)-\lambda^T(t_f)\mathbf{x}(t_f)]\nonumber
\\&&+\int_{t_0}^{t_f}[\mathcal{H}[\mathbf{x}(t),\mathbf{u}(t),\lambda(t),t]+\dot{\lambda}^T(t)\mathbf{x}(t)]dt\nonumber
\end{eqnarray}

After using the Hamiltonian, the stationarity of the cost function can be defined as variation of the cost function, $\Delta J$, being zero by any first order effect of control variation,$\Delta u$. This is similar in taking a Taylor series to the first order terms of J, shown in equation \eqref{71}

\begin{eqnarray}
\Delta J = \frac{\partial J}{\partial{\mathbf{u}}}\Delta \mathbf{u}+\frac{\partial J}{\partial{\mathbf{x}}}\Delta \mathbf{x}(\Delta \mathbf{u}) \label{71}
\end{eqnarray}

With a fixed final time, the first variation due to control is

\begin{eqnarray}
\Delta J = 
\bigg[\frac{\partial \phi}{\partial{\mathbf{x}}}-
\lambda^T]\Delta \mathbf{x}(\Delta \mathbf{u}) \bigg|_{t=t_f}+\lambda^T\Delta \mathbf{x}(\Delta \mathbf{u}) \bigg|_{t=t_0} +\int_{t_0}^{t_f}[\frac{\partial \mathcal{H}}{\partial{\mathbf{u}}}\Delta \mathbf{u}+[\frac{\partial \mathcal{H}}{\partial{\mathbf{x}}}+\dot{\lambda}^T]\Delta \mathbf{x}(\Delta \mathbf{u})\bigg]dt 
\end{eqnarray}

\begin{eqnarray}
\Delta J = \Delta J(t_f) + \Delta J(t_0) + \Delta J(t0,tf)
\end{eqnarray}

These 3 parts must equal zero separately for the control to be optimal. Luckily, we can assume that the control has no effect on the initial condition of the cost, so $\Delta J(t_0) = 0$. This leads to the next two equations that provide the basis for the necessary conditions of local optimality.

$\lambda(t)$ should satisfy the differential equation with the additional final condition.

\begin{eqnarray}
\dot{\lambda}^T(t) &=& -\frac{\partial \mathcal{H}}{\partial{\mathbf{x}}} \label{10}
\\\lambda^T(t_f) &=& \frac{\partial \phi}{\partial{\mathbf{x}}} \bigg|_{t=t_f} \label{11}
\end{eqnarray} 

Now that $\lambda$ has been identified, the next equation gives the equation to solve for the control function

\begin{eqnarray}
\frac{\partial \mathcal{H}}{\partial{\mathbf{u}}} = 0 \label{12}
\end{eqnarray}

Equations \eqref{10}, \eqref{11}, and \eqref{12} are called the Euler-Lagrange Equations. These are the necessary conditions when the Hamiltonian is stationary along the optimal trajectory. They are local conditions because they do not imply that there may be other optimal paths.

Because these equations are coupled, we may require a simpler method of solving them separately. One approach is to relate the Euler Lagrange Equations to the original cost function and the dynamics. $F$ and $G$ are the Jacobian matrices, $\frac{\partial{\mathbf{f}}}{\partial{\mathbf{x}}}$, and $\frac{\partial{\mathbf{f}}}{\partial{\mathbf{u}}}$.

\begin{eqnarray}
\frac{\partial \mathcal{H}}{\partial{\mathbf{x}}}^T = \frac{\partial \mathcal{L}}{\partial{\mathbf{x}}}^T + F^T\lambda
\\\frac{\partial \mathcal{H}}{\partial{\mathbf{u}}}^T = \frac{\partial \mathcal{L}}{\partial{\mathbf{u}}}^T + G^T\lambda
\end{eqnarray} 

The adjoint vector can be obtained by substituting into the dynamic equation with the final condition.

\begin{eqnarray}
\dot{\lambda}(t) = -F^T(t)\lambda(t)-\frac{\partial \mathcal{L}}{\partial{\mathbf{x}}}^T
\\ \lambda^T(t_f) = \frac{\partial \phi}{\partial{\mathbf{x}}} \bigg|_{t=t_f}
\end{eqnarray}

Finally, the optimal control history $\mathbf{u}$ can be calculated by solving the equation

\begin{eqnarray}
\frac{\partial \mathcal{L}}{\partial{\mathbf{u}}}^T+G^T(t)\lambda(t) = 0
\end{eqnarray}

Once again, this process utilizes the state dynamics as a constraint, a cost function based off of $\phi$ and $\mathcal{L}$, and calculus of variations to determine an optimal control function that will keep the energy required to drive the system from the initial state to the final state to a minimum. These are necessary conditions for optimality, but not sufficient conditions as the stationarity of the Hamiltonian may not necessarily be an optimal path.

In order to guarantee that the solution is a minimum for the optimal control problem, sufficient conditions  must be introduced. This section is primarily based off of the first of the three tests discussed in \cite{Stengel}. The test is the strengthened Legendre-Clebsch condition, which is shown below in equation \eqref{18}.

\begin{eqnarray}
\frac{\partial^2 \mathcal{H}}{\partial{\mathbf{u}}^2}[\mathbf{x^*},\mathbf{u^*},\mathbf{\lambda^*},t] = \mathcal{H}_u{_u} > 0 \label{18}
\end{eqnarray}

This equation can be interpreted that the Hamiltonian function is concave with respect to the control function. As a result, the stationary point that represents the zero variation of the Hamiltonian with respect to the control function is a minimum point, assuring optimality of the control function. This matrix is the Hessian matrix of the Hamiltonian with respect to the control. This sufficient condition is again only a local condition for the respective time period, globally the Hamiltonian may not always be concave.

\subsection{Simulated Systems}
Both of these optimal control frameworks are implemented on three canonical test systems for a variety of conditions. The first system these frameworks are simulated on is a cart-pole with 4 states and a single control. This under-actuated simulation a popular benchmark for trajectory optimization algorithms. The second test system is the double cart pole, which is a system with 6 states and still a single control. The ability for an optimal control algorithm to be able to optimize such a system is quite powerful. The final system is a quadrotor with 12 states and 4 controls. This is the beginning of a series of simulations that will eventually lead to a real-time implementation on a flying quadrotor.

\section{Differential Dynamic Programming}

Now we can apply the Dynamic Programming Principle to an algorithm for Differential Dynamic Programming. This algorithm is part of class of methods known as "shooting" direct methods. This is where the optimality conditions are not directly computed, and many trajectories are tested in order to find a minimum cost. Over each iteration, the trajectories are improved, leading to better cost.

\subsection{Equations}
Next we derive the basic equations required for the backward pass in Differential Dynamic Programming.

We first begin with the cost function that must be minimized $V(\mathbf{x}(t_0),t_0)$.

\begin{eqnarray}
V(\mathbf{x}(t_0),t_0) = \underset{\mathbf{u}}{\text{min}} \bigg[ \explain{\phi(\mathbf{x}(t_f),t_f)}{Cost to go} + \explain{ \int_{t_0}^{t_f}\mathcal{L}[\mathbf{x}(t),\mathbf{u}(t),t]dt}{Running Cost} \bigg]
\end{eqnarray}

This cost is subject to the dynamics of the system given an initial condition or a nominal trajectory ($\mathbf{\bar{x},\bar{u}}$, the algorithm of DDP must have a baseline to start from).

\begin{eqnarray}
\mathbf{\dot{x}}(t) = \mathbf{F}[\mathbf{x}(t),\mathbf{u}(t),t]
\end{eqnarray}

From my understanding, now we want to linearize the dynamics in order to produce an equation that will show the change with respect to a $\delta\mathbf{x}$ in terms of discrete time. This is required to generate a programmable algorithm that works in discrete time steps.

\begin{eqnarray}
\mathbf{\dot{x}}(t) =  \mathbf{F}[\mathbf{x}+\delta\mathbf{x},\mathbf{u}+\delta\mathbf{u},t] 
= \mathbf{F}[\mathbf{\bar{x}},\mathbf{\bar{u}},t] 
+ \nabla_x \mathbf{F} \delta\mathbf{x} + \nabla_u \mathbf{F} \delta\mathbf{u} 
\end{eqnarray}
\begin{eqnarray}
\explain{\mathbf{\dot{\delta x}} =
\mathbf{\dot{x}}(t) - \mathbf{F}[\mathbf{\bar{x}},\mathbf{\bar{u}},t]
= \nabla_x \mathbf{F} \delta\mathbf{x} + \nabla_u \mathbf{F} \delta\mathbf{u}}{Differential Dynamics}
\\\nonumber
\\\nonumber
\mathbf{\dot{\delta x}} = \frac{\mathbf{\delta x}(t_{k+1})-\mathbf{\delta x}(t_{k})}{dt} = \nabla_x \mathbf{F} \delta\mathbf{x}(t_k) + \nabla_u \mathbf{F}\delta\mathbf{u}(t_k)
\\\nonumber
\\\nonumber
\mathbf{\delta x}(t_{k+1}) = \mathbf{\delta x}(t_{k}) + \nabla_x \mathbf{F} \delta\mathbf{x}(t_k)dt + \nabla_u \mathbf{F}\delta\mathbf{u}(t_k)dt
\\\nonumber
\\\nonumber
\mathbf{\delta x}(t_{k+1}) 
= \explain{(I_{nxn}+\nabla_x \mathbf{F}dt)}{$\mathbf{\Phi}(t_k)$} \mathbf{\delta x}(t_{k})
+ \explain{(\nabla_u \mathbf{F}dt)}{$\mathbf{B}(t_k)$}\delta\mathbf{u}(t_k)
\end{eqnarray}

The completed linearized dynamics \eqref{lindam} along with the Bellman Dynamic Programming \eqref{bellman} Principle in discrete time gives the basis for Differential Dynamic Programming.

\begin{eqnarray}
\mathbf{\delta x}(t_{k+1})  = \mathbf{\Phi}(t_k) \mathbf{\delta x}(t_{k}) + \mathbf{B}(t_k) \delta\mathbf{u}(t_k) \label{lindam}
\\\nonumber
\\
V[\mathbf{x}(t_k),t_k]=\underset{\mathbf{u}}{\text{min}} \bigg[ \mathcal{L}[\mathbf{x}(t_k),\mathbf{u}(t_k),t_k]dt+V[\mathbf{x}(t_{k+1}),t_{k+1}]\bigg] \label{bellman}
\end{eqnarray} 

In our algorithm, the backwards pass is the one utilized to solve for the optimal control input. The backwards pass implies we must start from the desired state, and because Equation \eqref{bellman} contains that information in the form of $V[\mathbf{x}(t_{k+1}),t_{k+1}]$, the next step is to expand the value function $V$ out to the second order. By expanding this equation out, we can create a quadratic form for our value function with equations for $V$, $V_x$, and $V_{xx}$. After doing so we can solve for our optimal control.

\subsubsection*{Expand $V(\mathbf{x},t_{k+1})$ to the second order around the nominal trajectory}

\begin{eqnarray}
V(\mathbf{\bar{x}}(t_{k+1},t_{k+1})) =&& V \bigg( \mathbf{x}(t_{k+1}) - \mathbf{\bar{x}}(t_{k+1}
+ \mathbf{\bar{x}}(t_{k+1}) , t_{k+1} \bigg) \nonumber
\\\nonumber 
 =&& V \bigg( \mathbf{\bar{x}}(t_{k+1}) + \delta \mathbf{x}(t_{k+1}) ,t_{k+1} \bigg)
 \\\nonumber 
 =&& V \bigg( \mathbf{\bar{x}}(t_{k+1}), t_{k+1} \bigg) +
 \\
  && V_\mathbf{x}(t_{k+1})^T\delta\mathbf{x}(t_{k+1})
 + \frac{1}{2} \delta\mathbf{x}(t_{k+1})^T V_{\mathbf{xx}} (t_{k+1}) \delta\mathbf{x}(t_{k+1})
 \label{Vexp}
 \end{eqnarray}

In Equation \eqref{Vexp}, the following three terms are true.
\begin{eqnarray}
\delta \mathbf{x}(t_{k+1}) = \mathbf{x}(t_{k+1}) - \mathbf{\bar{x}}(t_{k+1})
\\
V_\mathbf{x}(t_{k+1}) = \nabla_\mathbf{x}V(\mathbf{x},t_{k+1}) \bigg|_{\mathbf{x} = \mathbf{\bar{x}}_{t_{k+1}}}
\\
V_\mathbf{xx}(t_{k+1}) = \nabla_\mathbf{xx}V(\mathbf{x},t_{k+1}) \bigg|_{\mathbf{x} = \mathbf{\bar{x}}_{t_{k+1}}}
\end{eqnarray}

Now we will substitute our equation for $\delta \mathbf{x}$, Equation \eqref{lindam} into the expansion for $V$, Equation \eqref{Vexp}. All values of $\mathbf{x}$ in the expansion are for the nominal value $\mathbf{\bar{x}}$.

\begin{eqnarray}
V(\mathbf{\bar{x}}(t_{k+1},t_{k+1})) =&& \nonumber 
V(t_{k+1})
+ V_\mathbf{x}(t_{k+1})^T \bigg( \mathbf{\Phi}(t_k) \mathbf{\delta x}(t_{k}) + \mathbf{B}(t_k) \delta\mathbf{u}(t_k) \bigg)+ 
\\\nonumber
 && \frac{1}{2} \bigg( \mathbf{\Phi}(t_k) \mathbf{\delta x}(t_{k}) + \mathbf{B}(t_k) \delta\mathbf{u}(t_k) \bigg)^T V_{\mathbf{xx}} (t_{k+1}) \bigg( \mathbf{\Phi}(t_k) \mathbf{\delta x}(t_{k}) + \mathbf{B}(t_k) \delta\mathbf{u}(t_k) \bigg)
 \\\nonumber
 \\\nonumber
 =&& V(t_{k+1}) + V_\mathbf{x}(t_{k+1})^T\mathbf{\Phi}(t_k)\mathbf{\delta x}(t_{k}) + V_\mathbf{x}(t_{k+1})^T \mathbf{B}(t_k) \delta\mathbf{u}(t_k)+
 \\\nonumber
 &&\frac{1}{2}\mathbf{\delta x}(t_{k})^T\mathbf{\Phi}(t_k)^T V_\mathbf{xx}(t_{k+1}) \mathbf{\Phi}(t_k) \mathbf{\delta x}(t_{k}) +
 \\\nonumber
 &&\frac{1}{2}\mathbf{\delta u}(t_{k})^T\mathbf{B}(t_k)^T V_\mathbf{xx}(t_{k+1}) \mathbf{B}(t_k) \mathbf{\delta u}(t_{k}) + 
 \\\nonumber
 && \frac{1}{2}\mathbf{\delta u}(t_{k})^T\mathbf{B}(t_k)^T V_\mathbf{xx}(t_{k+1}) \mathbf{\Phi}(t_k) \mathbf{\delta x}(t_{k}) +
 \\
 && \frac{1}{2}\mathbf{\delta x}(t_{k})^T\mathbf{\Phi}(t_k)^T V_\mathbf{xx}(t_{k+1}) \mathbf{B}(t_k) \mathbf{\delta u}(t_{k}) \label{Vfullexp}
\end{eqnarray}

\subsubsection*{Expand $\mathcal{L}(\mathbf{x},\mathbf{u},t_k)$ to the second order around the nominal trajectory}

This local quadratic expansion of $\mathcal{L}$ is very similar to the previous expansion of $V$. In this case, we have to take the derivatives of both $\mathbf{x}$ and $\mathbf{u}$.
\begin{eqnarray}
\mathcal{L} \bigg( \mathbf{x}(t_k),\mathbf{u}(t_k),t_k \bigg)dt = && \nonumber
\mathcal{L}\bigg( \mathbf{\bar{x}}(t_k) + \delta \mathbf{x}(t_k),\mathbf{\bar{u}}(t_k) + \delta \mathbf{u}(t_k),t_k \bigg)dt
\\\nonumber
=&& \explain{\mathcal{L}(\mathbf{x}(t_k),\mathbf{u}(t_k),t_k)dt}{$q(t_k)$} + \explain{(\nabla_\mathbf{x}\mathcal{L}(\mathbf{x}(t_k),\mathbf{u}(t_k),t_k)dt)^T}{$\mathbf{q}(t_k)^T$}\delta \mathbf{x}(t_k) +
\\\nonumber
&&\explain{(\nabla_\mathbf{u}\mathcal{L}(\mathbf{x}(t_k),\mathbf{u}(t_k),t_k)dt)^T}{$\mathbf{r}(t_k)^T$}\delta \mathbf{u}(t_k) + 
\\\nonumber 
&&\frac{1}{2} \delta \mathbf{x} (t_k)^T \explain{\nabla_\mathbf{xx}\mathcal{L}(\mathbf{x}(t_k),\mathbf{u}(t_k),t_k)dt}{$\mathbf{Q}(t_k)$}\delta \mathbf{x}(t_k) +
\\\nonumber 
&&\frac{1}{2} \delta \mathbf{u} (t_k)^T \explain{\nabla_\mathbf{uu}\mathcal{L}(\mathbf{x}(t_k),\mathbf{u}(t_k),t_k)dt}{$\mathbf{R}(t_k)$}\delta \mathbf{u}(t_k) +
\\\nonumber
&&\frac{1}{2} \delta \mathbf{u} (t_k)^T \explain{\nabla_\mathbf{xu}\mathcal{L}(\mathbf{x}(t_k),\mathbf{u}(t_k),t_k)dt}{$\mathbf{N}(t_k)$}\delta \mathbf{x}(t_k) +
\\
&&\frac{1}{2} \delta \mathbf{x} (t_k)^T \explain{\nabla_\mathbf{ux}\mathcal{L}(\mathbf{x}(t_k),\mathbf{u}(t_k),t_k)dt}{$\mathbf{M}(t_k)$}\delta \mathbf{u}(t_k)
\end{eqnarray}

\begin{eqnarray}
\mathcal{L} \bigg( \mathbf{x}(t_k),\mathbf{u}(t_k),t_k \bigg)dt =&& \nonumber
q(t_k) + \mathbf{q}(t_k)^T \delta \mathbf{x}(t_k) +
\mathbf{r}(t_k)^T \delta \mathbf{u}(t_k) +
\\\nonumber 
&& \frac{1}{2} \delta \mathbf{x} (t_k)^T \mathbf{Q}(t_k) \delta \mathbf{x}(t_k) + \frac{1}{2} \delta \mathbf{u} (t_k)^T \mathbf{R}(t_k) \delta \mathbf{u}(t_k) +
\\
&&\frac{1}{2} \delta \mathbf{u} (t_k)^T \mathbf{N}(t_k) \delta \mathbf{x}(t_k) +
\frac{1}{2} \delta \mathbf{x} (t_k)^T \mathbf{M}(t_k) \delta \mathbf{u}(t_k) \label{Lexp}
\end{eqnarray}

\subsubsection*{Combine the expansions of $\mathcal{L}$ and $V$}

Now we can combine Equations \eqref{Vfullexp} and \eqref{Lexp} in order to simplify our equations for the value function $V$.
\begin{eqnarray}
V[\mathbf{x}(t_k),t_k] =&&\underset{\delta\mathbf{u}(t_k)}{\text{min}} \bigg[ \explain{\mathcal{L}[\mathbf{x}(t_k),\mathbf{u}(t_k),t_k]dt}{Equation \eqref{Lexp}}+\explain{V[\mathbf{x}(t_{k+1}),t_{k+1}]}{Equation \eqref{Vfullexp}}\bigg] \nonumber
\\\nonumber
=&& \underset{\delta\mathbf{u}(t_k)}{\text{min}} \bigg[ q(t_k) + \mathbf{q}(t_k)^T \delta \mathbf{x}(t_k) +
\mathbf{r}(t_k)^T \delta \mathbf{u}(t_k) + \frac{1}{2} \delta \mathbf{x} (t_k)^T \mathbf{Q}(t_k) \delta \mathbf{x}(t_k) + 
\\\nonumber
&&\frac{1}{2} \delta \mathbf{u} (t_k)^T \mathbf{R}(t_k) \delta \mathbf{u}(t_k) +
\frac{1}{2} \delta \mathbf{u} (t_k)^T \mathbf{N}(t_k) \delta \mathbf{x}(t_k) +
\frac{1}{2} \delta \mathbf{x} (t_k)^T \mathbf{M}(t_k) \delta \mathbf{u}(t_k) +
\\\nonumber
&& V(t_{k+1}) + V_\mathbf{x}(t_{k+1})^T\mathbf{\Phi}(t_k)\mathbf{\delta x}(t_{k}) + V_\mathbf{x}(t_{k+1})^T \mathbf{B}(t_k) \delta\mathbf{u}(t_k)+
\\\nonumber
&&\frac{1}{2}\mathbf{\delta x}(t_{k})^T\mathbf{\Phi}(t_k)^T V_\mathbf{xx}(t_{k+1}) \mathbf{\Phi}(t_k) \mathbf{\delta x}(t_{k}) + \frac{1}{2}\mathbf{\delta u}(t_{k})^T\mathbf{B}(t_k)^T V_\mathbf{xx}(t_{k+1}) \mathbf{B}(t_k) \mathbf{\delta u}(t_{k}) + 
\\
&& \mathbf{\delta u}(t_{k})^T\mathbf{B}(t_k)^T V_\mathbf{xx}(t_{k+1}) \mathbf{\Phi}(t_k) \mathbf{\delta x}(t_{k}) \bigg] \label{FullV}
\end{eqnarray}

Now we can take Equation \eqref{FullV} and set the isolate the terms that contain $\delta \mathbf{u}(t_k)^T$. This will require some rearrangement of some terms such as $\mathbf{r}(t_k)$, but in the end, we will be able to solve for $\delta \mathbf{u}(t_k)$ because all the terms containing it must sum to zero.

\begin{eqnarray}
0 =&& \nonumber
\mathbf{r}(t_k) + \frac{1}{2}\mathbf{R}(t_k) \delta \mathbf{u}(t_k) + \frac{1}{2} \mathbf{N}(t_k) \delta \mathbf{x}(t_k) + \frac{1}{2} \mathbf{M}(t_k) \delta \mathbf{u}(t_k) + \mathbf{B}(t_k)^T V_\mathbf{x}(t_{k+1}) +
\\\nonumber
&& \mathbf{B}(t_k)^T V_\mathbf{xx}(t_{k+1}) \mathbf{B}(t_k) \mathbf{\delta u}(t_{k}) + \mathbf{B}(t_k)^T V_\mathbf{xx}(t_{k+1}) \mathbf{\Phi}(t_k) \mathbf{\delta x}(t_{k})
\end{eqnarray}

Now lets rearrange to put like terms together.

\begin{eqnarray}
0 =&& \nonumber
\explain{\bigg( \mathbf{r}(t_k) + \mathbf{B}(t_k)^T V_\mathbf{x}(t_{k+1}) \bigg)}{$\mathbf{g}$} + \explain{\bigg( \frac{1}{2}\mathbf{N}(t_k) + \frac{1}{2} \mathbf{M}(t_k)^T + \mathbf{B}(t_k)^T V_\mathbf{xx}(t_{k+1}) \mathbf{\Phi}(t_k) \bigg)}{$\mathbf{G}$} \mathbf{\delta x}(t_{k}) 
\\\nonumber
&&+ \frac{1}{2}\explain{\bigg( \mathbf{R}(t_k) + \mathbf{B}(t_k)^T V_\mathbf{xx}(t_{k+1}) \mathbf{B}(t_k) \bigg)}{$\mathbf{H}$} \delta \mathbf{u}(t_k)
\end{eqnarray}

\begin{eqnarray}
0 = \mathbf{g} + \mathbf{G}\mathbf{\delta x}(t_{k}) + \frac{1}{2}\mathbf{H}\delta \mathbf{u}(t_k) \nonumber
\\\nonumber
\\\nonumber
-\frac{1}{2}\mathbf{H}\delta \mathbf{u}(t_k) = \mathbf{g} + \mathbf{G}\mathbf{\delta x}(t_{k})
\\\nonumber
\\\nonumber
\delta \mathbf{u}(t_k) = \explain{-2\mathbf{H}^{-1} \mathbf{g}}{$\mathbf{l}(t_k)$} + \explain{-2\mathbf{H}^{-1} \mathbf{G}}{$\mathbf{L}(t_k)$} \mathbf{\delta x}(t_{k})
\\
\delta \mathbf{u}(t_k) = \mathbf{l}(t_k) + \mathbf{L}(t_k)\mathbf{\delta x}(t_{k}) \label{delu}
\end{eqnarray}

Now we can take Equation \eqref{delu} and substitute  for $\delta \mathbf{u}(t_k)$ in our expansion of V in Equation \eqref{FullV}.

\begin{eqnarray}
V[\mathbf{x}(t_k),t_k] =&& \nonumber  q(t_k) + \mathbf{q}(t_k)^T \delta \mathbf{x}(t_k) +
\mathbf{r}(t_k)^T \bigg( \mathbf{l}(t_k) + \mathbf{L}(t_k)\mathbf{\delta x}(t_{k}) \bigg) + \frac{1}{2} \delta \mathbf{x} (t_k)^T \mathbf{Q}(t_k) \delta \mathbf{x}(t_k) + 
\\\nonumber
&&\frac{1}{2} \bigg( \mathbf{l}(t_k) + \mathbf{L}(t_k)\mathbf{\delta x}(t_{k}) \bigg)^T \mathbf{R}(t_k) \bigg( \mathbf{l}(t_k) + \mathbf{L}(t_k)\mathbf{\delta x}(t_{k}) \bigg) +
\\\nonumber
&&\frac{1}{2} \bigg( \mathbf{l}(t_k) + \mathbf{L}(t_k)\mathbf{\delta x}(t_{k}) \bigg)^T \mathbf{N}(t_k) \delta \mathbf{x}(t_k) +
\frac{1}{2} \delta \mathbf{x} (t_k)^T \mathbf{M}(t_k) \bigg( \mathbf{l}(t_k) + \mathbf{L}(t_k)\mathbf{\delta x}(t_{k}) \bigg) +
\\\nonumber
&& V(t_{k+1}) + V_\mathbf{x}(t_{k+1})^T\mathbf{\Phi}(t_k)\mathbf{\delta x}(t_{k}) + V_\mathbf{x}(t_{k+1})^T \mathbf{B}(t_k)\bigg( \mathbf{l}(t_k) + \mathbf{L}(t_k)\mathbf{\delta x}(t_{k}) \bigg)+
\\\nonumber
&&\frac{1}{2}\mathbf{\delta x}(t_{k})^T\mathbf{\Phi}(t_k)^T V_\mathbf{xx}(t_{k+1}) \mathbf{\Phi}(t_k) \mathbf{\delta x}(t_{k}) + 
\\\nonumber
&&\frac{1}{2}\bigg( \mathbf{l}(t_k) + \mathbf{L}(t_k)\mathbf{\delta x}(t_{k}) \bigg)^T\mathbf{B}(t_k)^T V_\mathbf{xx}(t_{k+1}) \mathbf{B}(t_k) \bigg( \mathbf{l}(t_k) + \mathbf{L}(t_k)\mathbf{\delta x}(t_{k}) \bigg) + 
\\
&& \bigg( \mathbf{l}(t_k) + \mathbf{L}(t_k)\mathbf{\delta x}(t_{k}) \bigg)^T\mathbf{B}(t_k)^T V_\mathbf{xx}(t_{k+1}) \mathbf{\Phi}(t_k) \mathbf{\delta x}(t_{k}) \label{finalform}
\end{eqnarray}

\subsubsection*{Split into Zero, First, and Second Order Groupings}

For the final step of this process, each term in Equation \eqref{finalform} can be separated into terms that are zero, first, and second order with $\delta \mathbf{x}$. This step allows us to derive the Ricatti Equations for $V(t_k)$, $V_\mathbf{x}(t_k)$, and $V_{\mathbf{xx}}(t_k)$. 

For each of the equations, $\delta \mathbf{x}$ is then set to zero. $\delta \mathbf{x} =0$ because $V$ signifies a minimum, and at the minimum the variation in value with respect to changes in state should be zero. We already reflected the variation with respect to the control when we solved for $\delta \mathbf{u}(t_k)$ in Equation \eqref{delu}.

\subsubsection{Zero Order Terms}

\begin{eqnarray}
=&& q(t_k) + \mathbf{r}(t_k)^T \mathbf{l}(t_k) + \frac{1}{2} \mathbf{l}(t_k)^T\mathbf{R}(t_k)\mathbf{l}(t_k) + V(t_{k+1}) + V_\mathbf{x}(t_{k+1})^T \mathbf{B}(t_k) \mathbf{l}(t_k)
\\\nonumber
&& + \frac{1}{2} \mathbf{l}(t_k)^T\mathbf{B}(t_k)^T V_\mathbf{xx}(t_{k+1}) \mathbf{B}(t_k)\mathbf{l}(t_k)
\\\nonumber
=&& \nonumber
q(t_k) +  V(t_{k+1}) + \explain{\bigg( \mathbf{r}(t_k)^T + V_\mathbf{x}(t_{k+1})^T \mathbf{B}(t_k) \bigg)}{$\mathbf{g}(t_k)^T$}\mathbf{l}(t_k) +
\\\nonumber
&&\frac{1}{2}\mathbf{l}(t_k)^T\explain{\bigg( \mathbf{R}(t_k) + \mathbf{B}(t_k)^T V_\mathbf{xx}(t_{k+1}) \mathbf{B}(t_k) \bigg)}{$\mathbf{H}(t_k)$}\mathbf{l}(t_k)
\\
=&& q(t_k) +  V(t_{k+1}) + \mathbf{g}(t_k)^T\mathbf{l}(t_k) + \frac{1}{2}\mathbf{l}(t_k)^T\mathbf{H}(t_k)\mathbf{l}(t_k) \label{zero}
\end{eqnarray}

\subsubsection*{First Order Terms}
\begin{eqnarray}
=&& \nonumber \mathbf{q}(t_k)^T\delta \mathbf{x}(t_k) + \mathbf{r}(t_k)^T\mathbf{L}(t_k)\delta \mathbf{x}(t_k) + \mathbf{l}(t_k)^T\mathbf{R}(t_k)\mathbf{L}(t_k)\delta \mathbf{x}(t_k) + \frac{1}{2}\mathbf{l}(t_k)^T\mathbf{N}(t_k)\delta \mathbf{x}(t_k)+
\\\nonumber
\\\nonumber
&& \frac{1}{2}\mathbf{l}(t_k)\mathbf{M}(t_k)^T\delta \mathbf{x}(t_k) +V_\mathbf{x}(t_{k+1})^T\mathbf{\Phi}(t_k)\mathbf{\delta x}(t_{k}) + V_\mathbf{x}(t_{k+1})^T \mathbf{B}(t_k)\mathbf{L}(t_k)\mathbf{\delta x}(t_{k})+
\\\nonumber
\\
&& \mathbf{l}(t_k)^T\mathbf{B}(t_k)^T V_\mathbf{xx}(t_{k+1}) \mathbf{B}(t_k)\mathbf{L}(t_k)\mathbf{\delta x}(t_{k}) + \mathbf{l}(t_k)^T\mathbf{B}(t_k)^TV_\mathbf{xx}(t_{k+1}) \mathbf{\Phi}(t_k) \mathbf{\delta x}(t_{k})
\\\nonumber
\\\nonumber
=&& \bigg(\mathbf{q}(t_k)^T + V_\mathbf{x}(t_{k+1})^T\mathbf{\Phi}(t_k)  \bigg)\delta\mathbf{x}(t_k) + \explain{\bigg( \mathbf{r}(t_k)^T + V_\mathbf{x}(t_{k+1})^T \mathbf{B}(t_k) \bigg)}{$\mathbf{g}(t_k)^T$}\mathbf{L}(t_k)\mathbf{\delta x}(t_{k}) +
\\\nonumber
&&\mathbf{l}(t_k)^T \explain{\bigg( \mathbf{R}(t_k) + \mathbf{B}(t_k)^T V_\mathbf{xx}(t_{k+1}) \mathbf{B}(t_k) \bigg)}{$\mathbf{H}(t_k)$}\mathbf{L}(t_k)\delta \mathbf{x}(t_k) + 
\\\nonumber
&&\mathbf{l}(t_k)^T \explain{\bigg( \frac{1}{2}\mathbf{N}(t_k) + \frac{1}{2}\mathbf{M}(t_k)^T +\mathbf{B}(t_k)^TV_\mathbf{xx}(t_{k+1}) \mathbf{\Phi}(t_k) \bigg)}{$\mathbf{G}(t_k)$}\delta \mathbf{x}(t_k)
\\\nonumber
\\\nonumber
=&&\bigg(\mathbf{q}(t_k)^T + V_\mathbf{x}(t_{k+1})^T\mathbf{\Phi}(t_k)  \bigg)\delta\mathbf{x}(t_k) + \mathbf{g}(t_k)^T\mathbf{L}(t_k)\mathbf{\delta x}(t_{k})+
\\\nonumber
\\
&& \mathbf{l}(t_k)^T\mathbf{H}(t_k)\mathbf{L}(t_k)\delta \mathbf{x}(t_k) + \mathbf{l}(t_k)^T\mathbf{G}(t_k)\delta \mathbf{x}(t_k)
\end{eqnarray}

\subsubsection*{Second Order Terms}
\begin{eqnarray}
=&& \nonumber \frac{1}{2} \delta \mathbf{x} (t_k)^T \mathbf{Q}(t_k) \delta \mathbf{x}(t_k) +\frac{1}{2}\mathbf{\delta x}(t_{k})^T\mathbf{L}(t_k)^T\mathbf{R}(t_k)\mathbf{L}(t_k)\delta \mathbf{x}(t_k) +
\\\nonumber
\\\nonumber
&& \frac{1}{2}\mathbf{\delta x}(t_{k})^T\mathbf{L}(t_k)^T\mathbf{N}(t_k)\delta \mathbf{x}(t_k) +\frac{1}{2}\mathbf{\delta x}(t_{k})^T\mathbf{M}(t_k)\mathbf{L}(t_k)\delta\mathbf{x}(t_k)+
\\\nonumber
\\\nonumber
&& \frac{1}{2}\mathbf{\delta x}(t_{k})^T\mathbf{\Phi}(t_k)^T V_\mathbf{xx}(t_{k+1}) \mathbf{\Phi}(t_k) \mathbf{\delta x}(t_{k})+
\\\nonumber
\\\nonumber
&& \frac{1}{2}\mathbf{\delta x}(t_{k})^T\mathbf{L}(t_k)^T\mathbf{B}(t_k)^TV_\mathbf{xx}(t_{k+1})\mathbf{B}(t_k)\mathbf{L}(t_k)\delta \mathbf{x}(t_k) +
\\\nonumber
\\
&&\mathbf{\delta x}(t_{k})^T\mathbf{L}(t_k)^T\mathbf{B}(t_k)^TV_\mathbf{xx}(t_{k+1})\mathbf{\Phi}(t_k)\mathbf{\delta x}(t_{k})
\\\nonumber
\\\nonumber
=&&\nonumber \frac{1}{2} \delta \mathbf{x} (t_k)^T \mathbf{Q}(t_k) \delta \mathbf{x}(t_k) +
\\\nonumber
\\\nonumber
&&\frac{1}{2}\mathbf{\delta x}(t_{k})^T\mathbf{L}(t_k)^T\explain{\bigg( \mathbf{R}(t_k) + \mathbf{B}(t_k)^T V_\mathbf{xx}(t_{k+1}) \mathbf{B}(t_k) \bigg)}{$\mathbf{H}(t_k)$}\mathbf{L}(t_k)\delta \mathbf{x}(t_k) +
\\\nonumber
\\\nonumber
&&\frac{1}{2}\mathbf{\delta x}(t_{k})^T\mathbf{L}(t_k)^T\bigg( \mathbf{N}(t_k) + \mathbf{B}(t_k)^TV_\mathbf{xx}(t_{k+1})\mathbf{\Phi}(t_k) \bigg)\delta \mathbf{x}(t_k)+
\\\nonumber
\\\nonumber
&&\frac{1}{2}\mathbf{\delta x}(t_{k})^T\bigg( \mathbf{M}(t_k) + \mathbf{\Phi}(t_k)V_\mathbf{xx}(t_{k+1})\mathbf{B}(t_k)^T \bigg)\mathbf{L}(t_k)\delta+ \mathbf{x}(t_k)
\\\nonumber
\\
&&\frac{1}{2}\mathbf{\delta x}(t_{k})^T\mathbf{\Phi}(t_k)^T V_\mathbf{xx}(t_{k+1}) \mathbf{\Phi}(t_k) \mathbf{\delta x}(t_{k}
\end{eqnarray}

\subsubsection*{Riccati Equations}
After the terms are separated by order, derivatives can be taken in $\delta \mathbf{x}$ in order to get the respective Riccati Equations, shown below.

\begin{eqnarray}
V_\mathbf{xx}(t_k)=&& \nonumber \mathbf{Q}(t_k) + \mathbf{L}(t_k)^T\mathbf{H}(t_k)\mathbf{L}(t_k) + \mathbf{L}(t_k)^T\bigg( \mathbf{N}(t_k) + \mathbf{B}(t_k)^TV_\mathbf{xx}(t_{k+1})\mathbf{\Phi}(t_k) \bigg)\mathbf{L}(t_k)
\\\nonumber
&&\mathbf{L}(t_k)^T\bigg( \mathbf{M}(t_k) + \mathbf{\Phi}(t_k)V_\mathbf{xx}(t_{k+1})\mathbf{B}(t_k)^T \bigg)\mathbf{L}(t_k) + \mathbf{\Phi}(t_k)^T V_\mathbf{xx}(t_{k+1}) \mathbf{\Phi}(t_k)
\\\nonumber
\\\nonumber
V_\mathbf{x}(t_k)=&& \mathbf{q}(t_k)^T + V_\mathbf{x}(t_{k+1})^T\mathbf{\Phi}(t_k) + \mathbf{g}(t_k)^T\mathbf{L}(t_k) +  \mathbf{l}(t_k)^T\mathbf{H}(t_k)\mathbf{L}(t_k) + \mathbf{l}(t_k)^T\mathbf{G}(t_k)
\\\nonumber
\\\nonumber
V(t_k) =&& q(t_k) +  V(t_{k+1}) + \mathbf{g}(t_k)^T\mathbf{l}(t_k) + \frac{1}{2}\mathbf{l}(t_k)^T\mathbf{H}(t_k)\mathbf{l}(t_k)
\end{eqnarray}

In the following section we will review the equations and methods for Gauss Pseudospectral Optimal Control, then compare the results of both algorithms.

\section{Gaussian Pseudospectral Optimal Control}
Pseudospectral optimal control methods are a significant deviation from the classical methods for optimal control. Typically, classical "indirect" methods are known to analytically solve for first and second order optimality conditions. Then a cost function in optimized by solving for a control that satisfies the dynamic constraints and yields the lowest cost. Unfortunately, find these analytical solutions and evaluating the optimality of each solution can be difficult.

On the other hand, Pseudospectral optimal control is in a class of algorithms known as "direct" methods. The general outline (from Huntington) is to do the following: First, We will convert the dynamic system into a problem with a finite set of variable and transform the the constraints to ensure they are simply algebraic constraints. This will allow us to solve with a Non-linear Program (NLP) Solver. Second, We will solve this finite dimensional problem with the NLP and optimize the parameters with which we approximated the original system. Third, we can assess the accuracy of our approximation utilizing the costate and the Karush-Kuhn-Tucker (KKT) conditions. The advantage of this particular method (for smooth problems) is that it exhibits a fast convergence rate. The Gauss Pseudospectral Method is described in much further detail in \cite{Huntington}.
\subsection{Equations}
Let us begin with a few equations to describe the problem. The Pseudospectral methods are defined with the Transformed Continuous Bolza problem.

\begin{eqnarray}
J = \Phi(\mathbf{x}(\tau_0),t_0,\mathbf{x}(\tau_f),t_f) + \frac{t_f-t_0}{2} \int_{\tau_0}^{\tau_f}g(\mathbf{x}(\tau),\mathbf{u}(\tau),\tau;t_0,t_f)d\tau \label{Cost}
\end{eqnarray}

This cost will be minimized subject to the following constraints.
\begin{equation}
\frac{d\mathbf{x}}{dt} = \frac{t_f-t_0}{2}\mathbf{f}(\mathbf{x}(\tau),\mathbf{u}(\tau),\tau;t_0,t_f) \label{dyn}
\end{equation}
\begin{equation}
\phi(\mathbf{x}(\tau_0),t_0,\mathbf{x}(\tau_f),t_f) = \mathbf{0} \label{bound}
\end{equation}
\begin{equation}
\mathbf{C}(\mathbf{x}(\tau),\mathbf{u}(\tau),\tau;t_0,t_f) \leq \mathbf{0} \label{path}
\end{equation}
\begin{eqnarray}
\tau = \frac{2t}{t_f-t_0}-\frac{t_f+t_0}{t_f-t_0} \label{time}
\end{eqnarray}

Equation (\ref{dyn}) is the dynamic constraint to the problem. Equation (\ref{bound}) is the boundary condition at the final state. Equation (\ref{path}) is the path constraint for the problem. Equation (\ref{time}) shows the time transformation required for this derivation of pseudospectral optimal control. Typically, we require a fixed time interval such as $[-1,1]$, thus such a transformation can be used that is still valid with free initial and final times $t_0$ and $t_f$.

\subsubsection*{Collocation Point Search}

One of the most important aspects of this method is the choice of collocation points $K$. These points are where we set the approximation equal to the function across the interval and provided the set for our approximation. While there are a variety of points to choose from, the three most accurate choices are Legendre-Gauss (LG), Legendre-Gauss-Radau (LGR), and Legendre-Gauss-Lobatto (LGL) Points. Below is the equation for the Legendre-Gauss points. The points are defined as the roots of Equation (\ref{LG}). Equation (\ref{LG}) is the $K^{th}$ degree Legendre polynomial.

\begin{equation}
P_K(\tau) = \frac{1}{2^KK!}\frac{d^K}{d\tau^K}\big[(\tau^2-1)^K \big] \label{LG}
\end{equation}

The reason for the points to be appropriated as the roots of this polynomial is because of accuracy of the approximation. Equally spaced points will not accurately represent the function near the boundaries, hence this distribution actually clusters the points closer to the boundaries of $[-1,1]$. This is another reason of why we need to have a finite boundary for this method (in the transformed space).

\subsubsection*{State, Control, and Costate Approximation}

The next step of this is to define how we will approximate the state, control, and the costate. Note that we do not have to use the same approximation for each one (as in, we do not have to use the same collocation points and basis functions for each state, control, and costate), but we do use the same process. From Huntington, "Let us assume there are $K$ collocation points" (as defined above) and use these points to estimate both the state and control.

\begin{equation}
x(\tau) \approx \mathbf{X}(\tau) = \sum_{i=1}^K \mathcal{L}_i(\tau)\mathbf{X}(\tau_i) \label{state}
\end{equation}

\begin{equation}
\mathcal{L}_i(\tau) = \prod_{j=1,j\neq i}^{K} \frac{\tau-\tau_j}{\tau_i-\tau_j} = \frac{g(\tau)}{(\tau-\tau_i) \dot{g}}(\tau) \label{Lag}
\end{equation}

\begin{equation}
g(\tau) = (1+\tau)P_k(\tau)
\end{equation}

Equation (\ref{state}) represents the state approximation in terms of a set of Lagrange interpolating polynomials. Equation (\ref{Lag}) represents the actual polynomials, where the $P_K$ is the Legendre Polynomial. 

\subsubsection*{Integral Approximation via Gauss Quadrature}

Within the dynamic constraints and cost function we also have integrals that must be approximated. One way to do this is via a Gauss quadrature. The generic form for this is in Equation (\ref{GQ}).

\begin{equation}
\int_{a}^{b}f(\tau)d\tau \approx \sum_{i=1}^{K}w_if(\tau_i) \label{GQ}
\end{equation}

The points $\tau_k$ are the quadrature points on the interval $[-1,1]$ and the weights $w_i$ are the quadrature weights. Equation (\ref{qweight}) are for finding these weights.

\begin{equation}
w_i = \int_{-1}^{1}\mathcal{L}_i(\tau)d\tau = \frac{2}{(1-\tau_i^2)[\dot{P}_k(\tau_i)]^2},  (i = 1,...,K) \label{qweight}
\end{equation}

\subsubsection*{Algrebriac Representation of Dynamics}

The next step is to convert the differential dynamic equation constraint on the optimization problem into an algebraic constraint. This is done in Equation(\ref{DtoA})

\begin{equation}
\dot{x}(\tau_k) \approx \dot{\mathbf{X}}(\tau_k) = \sum_{i=1}^{K}\dot{\mathcal{L}}(\tau_k)\mathbf{X}(\tau_i) = \sum_{i=1}^{K}D_{ki}\mathbf{X}(\tau_k), (k = 1,...,K) \label{DtoA}
\end{equation}

The differentiation matrix which is given in Equation (\ref{Diff}).
\begin{equation}
D_{ki} = \begin{cases}
\frac{(1+\tau_k)\dot{P}_K(\tau_k)+P_K(\tau_k)}{(\tau_k-\tau_i)\big[(1+\tau_i)\dot{P}_K(\tau_i)+P_K(\tau_i) \big]} &\text{if } k\neq i\\
\\
\frac{(1+\tau_i)\ddot{P}_K(\tau_i)+2\dot{P}_K(\tau_i)}{2\big[ (1+\tau_i)\dot{P}_K(\tau_i)+P_K(\tau_i) \big]} & \text{if } k = i

\end{cases} \label{Diff}
\end{equation}

We can now implement this algebraic representation in terms of a residual function by equation the derivative of our approximation to the vector field.

\begin{equation}
\mathbf{R}_k = \sum_{i=1}^{K}D_{ki}\mathbf{X}(\tau_k) - \frac{t_f-t_0}{2}\mathbf{f}(\mathbf{X}(\tau_k),\mathbf{\mathbf{U}}(\tau_k),\tau_k;t_0,t_f) = \mathbf{0},  (k = 1,...,K) \label{res}
\end{equation}

\subsubsection*{Gauss Pseudospectral Discretization of the Continuous Bolza Problem}

Now we can put everything together and define our new problem. We want to minimize the following cost function Equation (\ref{DBP}) subject to the algebraic collocation constraints in Equation (\ref{AC}). We also have the quadrature constraint in Equation (\ref{QC}), boundary constraint in Equation (\ref{BC}) and path constraint in Equation (\ref{PC}).

\begin{equation}
J = \Phi(\mathbf{X}_0,t0,\mathbf{X}_f,t_f) + \frac{t_f-t_0}{2}\sum_{k=1}^{K}w_kg(\mathbf{X}_k,\mathbf{U}_k,\tau_k;t_0,t_f) \label{DBP}
\end{equation}

\begin{equation}
\mathbf{R}_k = \sum_{i=0}^{K}D_{ki}\mathbf{X}_i-\frac{t_f-t_0}{2}\mathbf{f}(\mathbf{X}_k,\mathbf{U}_k,\tau_k;t_0,t_f) = \mathbf{0}  (k = 1,...,N) \label{AC}
\end{equation}

\begin{equation}
\mathbf{R}_f = \mathbf{X}_f-\mathbf{X}_0-\frac{t_f-t_0}{2}\sum_{k=1}^{K}w_k\mathbf{f}(\mathbf{X}_k,\mathbf{U}_k,\tau_k;t_0,t_f) = \mathbf{0} \label{QC}
\end{equation}

\begin{equation}
\phi(\mathbf{X}_0,t_0,\mathbf{X}_f,t_f) = \mathbf{0} \label{BC}
\end{equation}

\begin{equation}
\mathbf{C}(\mathbf{X}_k,\mathbf{U}_k,\tau_k;t_0,t_f) \leq \mathbf{0} \label{PC}
\end{equation}

The procedure to solve is as follows.

\begin{enumerate}
\item Choose the initial and final times, then choose your basis functions and use the Legendre polynomials to get the collocation points.
\item Translate the dynamic system into the the discretized and approximated version of the Continuous Bolza Problem. Also translate the constraints into the equivalent forms.
\item Solve the finite-dimensional problem with a Non-linear Program Solver.
\item Check the accuracy of the problem by the costate approximation (not listed here), and then run again with different conditions, moving towards the optimum.
\end{enumerate}

The program utilized to implement this method is called GPOPS II \cite{GPOPS}.The non-linear program was solved using \textit{ipopt}.

\section{Results}
We applied both algorithms to three systems. The first was a simple cart pole, the second was a double cart pole, and the final system was a quadrotor. Table \ref{Problem} below contains the parameters and cost functions for each problem. The parameters and cost functions were kept consistent between both algorithms. We will observe the behavior of each algorithm, as well as the final cost and total runtime for each algorithm. Note that is not possible to directly compare the run time as DDP is currently implemented in MATLAB, while GPOC has a few files that run through C++. Nevertheless, a few generic observations can be made about algorithm speed. The cost comparison is only approximate, as DDP is a discrete time optimizer (i.e. the dynamics are discretized and linearized locally at each timestep) and GPOC is a continuous time optimizer. With a small timestep $dt$ for DDP, the costs are comparable.

\begin{table}[H]
\centering
\caption{Problem Settings}
\label{Problem}
\begin{tabular}{l | l | c | c}
Problem & System Parameters & Final Time & Cost Function\\ \hline \hline
\multirow{3}{5em}{Cart Pole} & Cart Mass = 1 kg  & &
						 \\ & Link Mass = 5 kg & 2 seconds & $\int_{t_0}^{tf}(Ru^2)dt $
					 	 \\ & Link Length = 1.5 m & & \\ \hline
\multirow{5}{5em}{Double Cart Pole} & Cart Mass = 3 kg & &
								\\ & Link 1 Mass = 1 kg & &
								\\ & Link 2 Mass = 20 kg & 4 seconds & $\int_{t_0}^{tf}(Ru^2)dt $
								\\ & Link 1 Length = 1.5 m & &
								\\ & Link 2 Length = 1.5 m & & \\ \hline
\multirow{6}{5em}{Quadrotor} & Quadrotor Mass = 1 kg & &
						 \\ & $J_x$ = 8.1E-3& &
						 \\ & $J_x$ = 8.1E-3 & 3 seconds & $\int_{t_0}^{tf}(x^TQx + u^TRu)dt $
						 \\ & $J_x$ = 14.2E-3 & &
						 \\ & Arm Length = .24 m
\end{tabular}
\end{table}

\subsection{Cart Pole}
\begin{figure}[H]
	\centering
	\begin{subfigure}[b]{0.48\textwidth}
		\includegraphics[width=\textwidth]{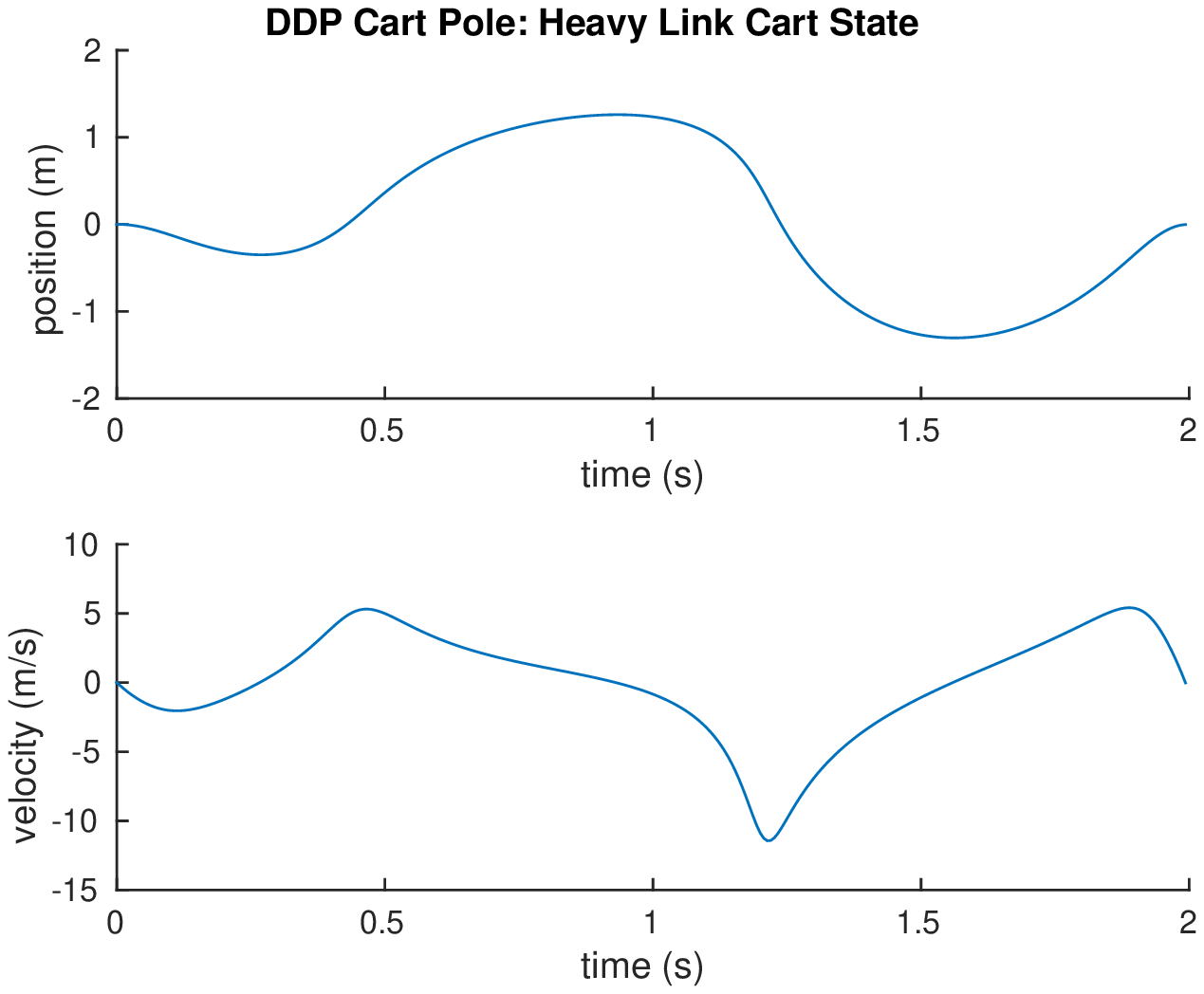}
		\caption{}
		\label{DCPCS}
	\end{subfigure}
	~
	\begin{subfigure}[b]{0.48\textwidth}
		\includegraphics[width=\textwidth]{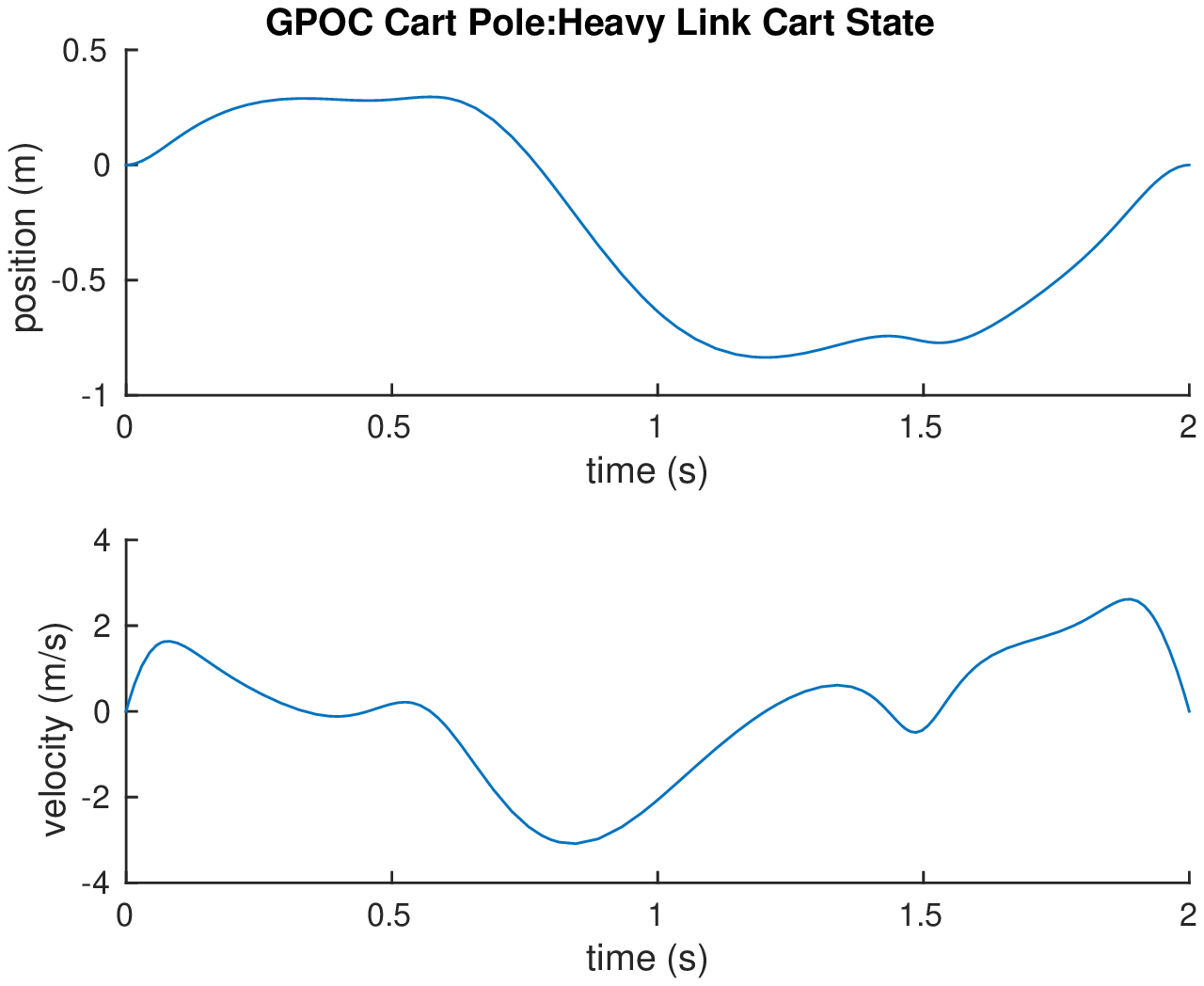}
		\caption{}
		\label{GCPCS}
	\end{subfigure}
	\\
		\begin{subfigure}[b]{0.48\textwidth}
		\includegraphics[width=\textwidth]{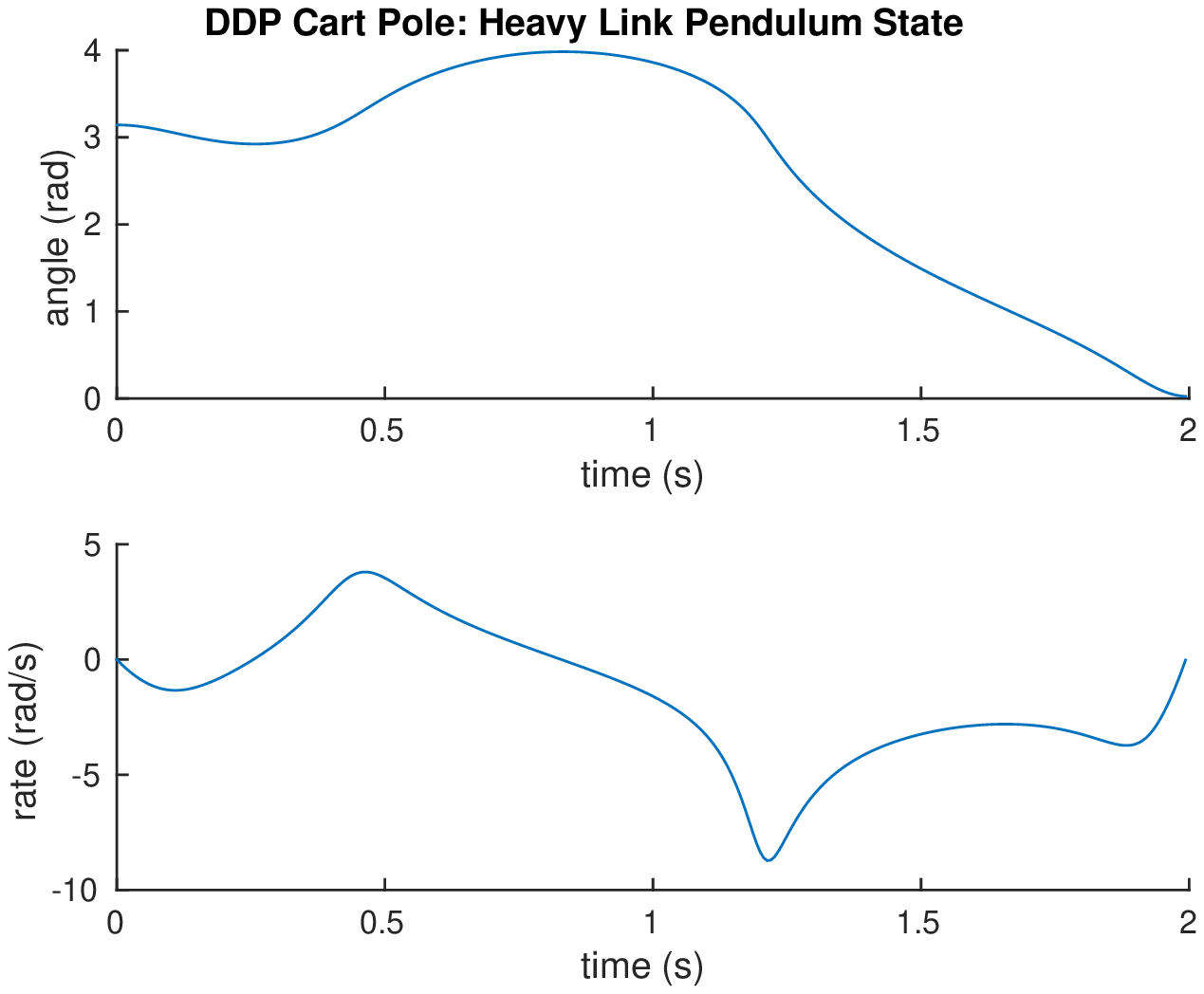}
		\caption{}
		\label{DCPPS}
	\end{subfigure}
	~
	\begin{subfigure}[b]{0.48\textwidth}
		\includegraphics[width=\textwidth]{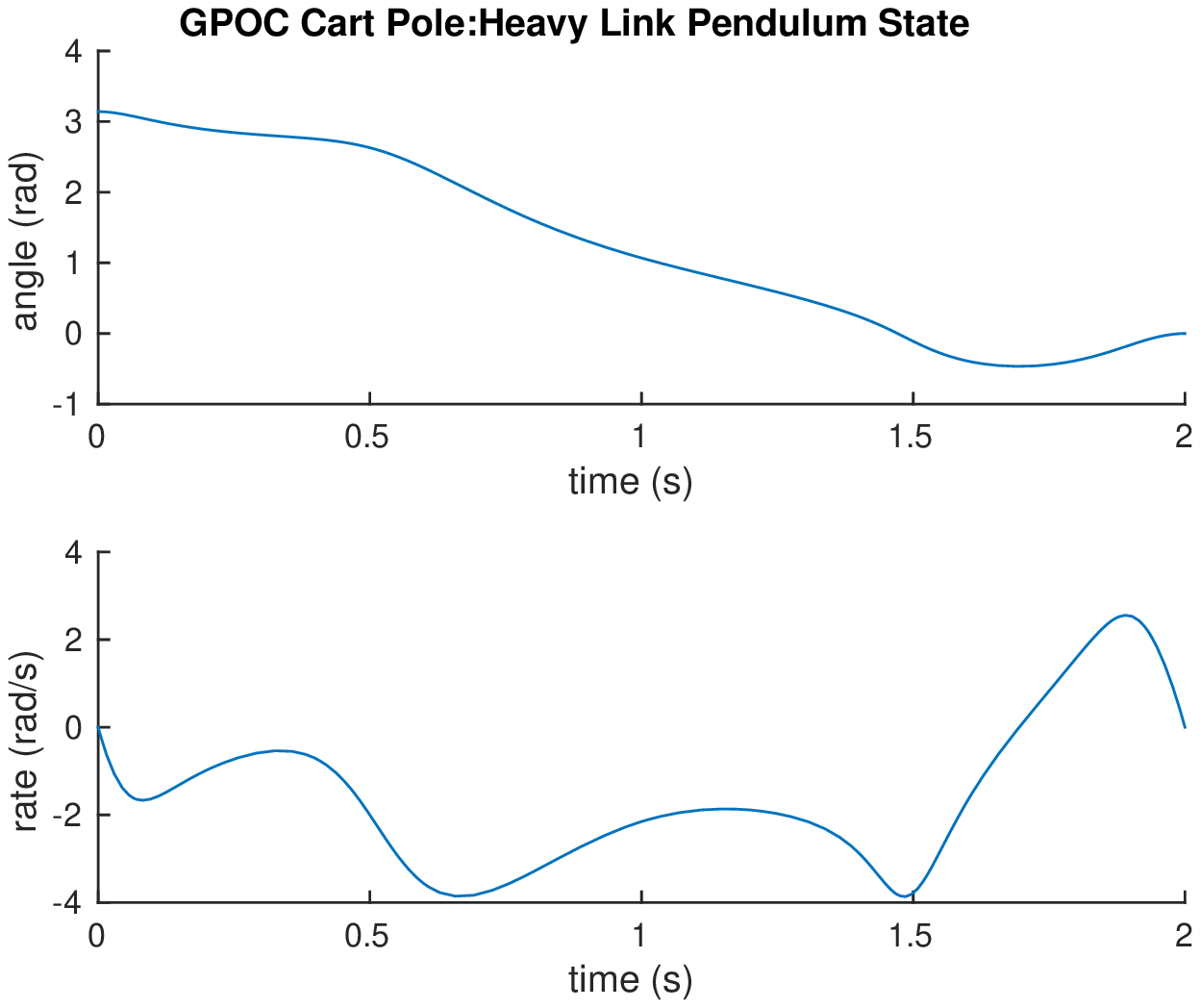}
		\caption{}
		\label{GCPPS}
	\end{subfigure}
	\caption{Cart Pole State Comparison}
	\label{CPCS}
	\begin{subfigure}[b]{0.48\textwidth}
		\includegraphics[width=\textwidth]{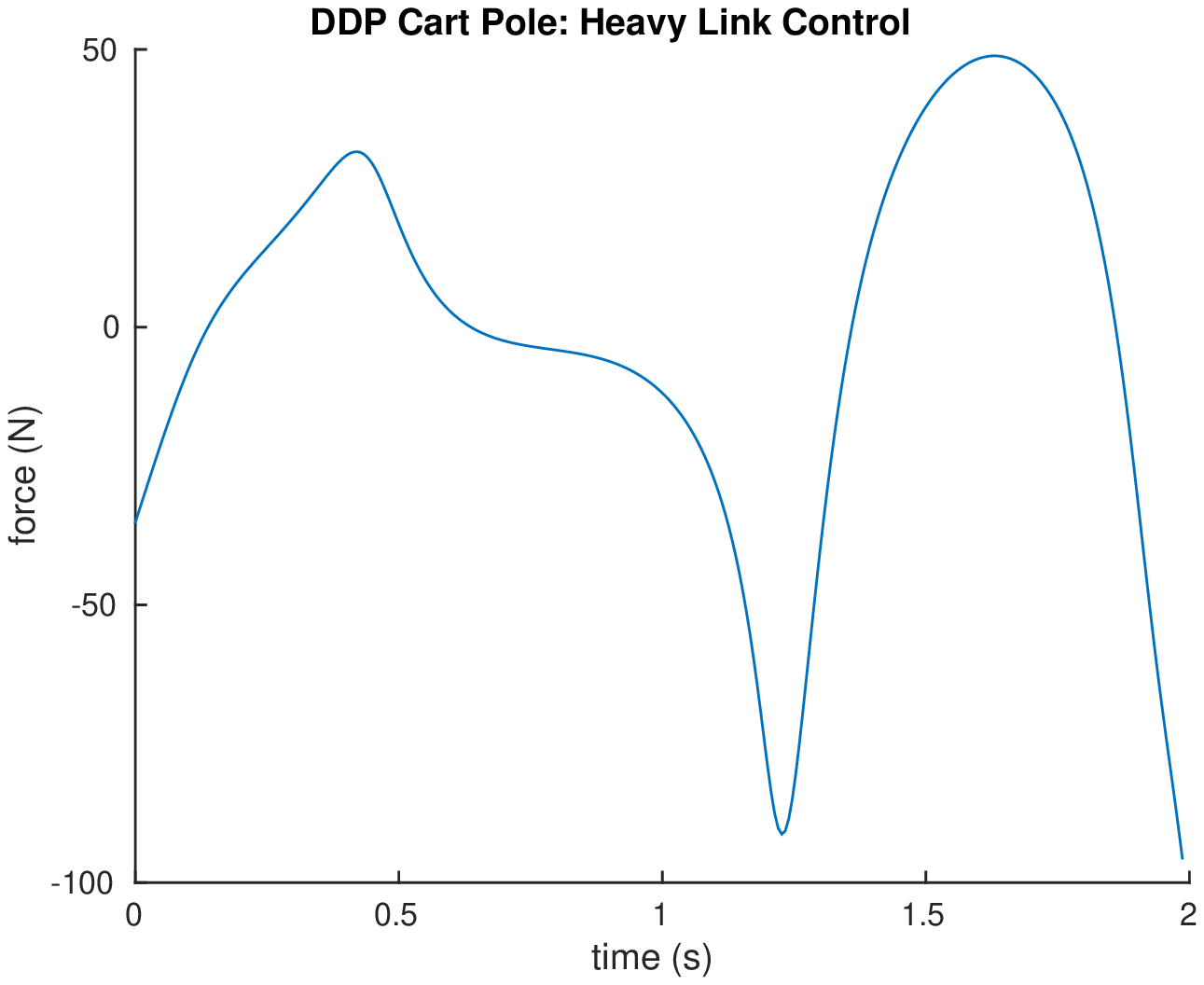}
		\caption{}
		\label{DCPControl}
	\end{subfigure}
	~
	\begin{subfigure}[b]{0.48\textwidth}
		\includegraphics[width=\textwidth]{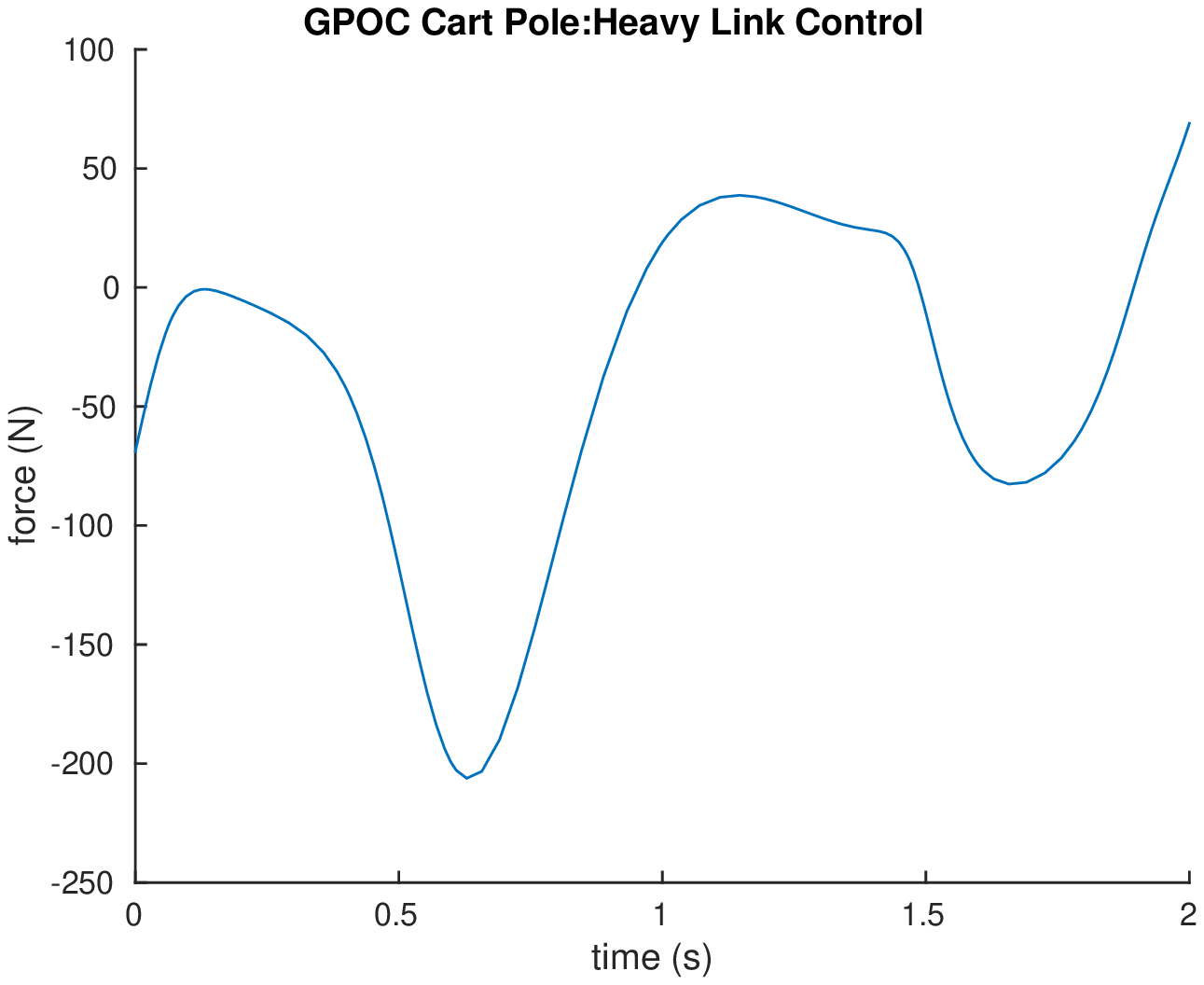}
		\caption{}
		\label{GCPControl}
	\end{subfigure}
	\caption{Cart Pole Control Comparison}
\end{figure}

The first system is the classic cart pole. The task is to swing up the cart pole from a downwards stable equilibrium to an unstable upwards equilibrium. In Figure \ref{CPCS}, the comparison between the Differential Dynamic Programming and Gauss Pseudospectral method can be observed. The left column corresponds to DDP, while the right column corresponds to GPOC. In this first comparison, DDP ended with a final cost of 0.756, while GPOC ended with a final cost of 0.1257. The runtime for DDP was 7.5 seconds, while GPOC ran in 0.9 seconds. In Figure \ref{DCPCS}, we can see that the cart moved to the left and right of the origin far more than for in Figure \ref{GCPCS}. This behavior is characteristic of GPOC, it often decides that a higher control effort is more optimal in order to reach the target. The observation of the angular rate bounds can also be observed is Figures \ref{DCPPS} and \ref{GCPPS}.

\begin{figure}[H]
\centering
	\includegraphics[width=.48\textwidth]{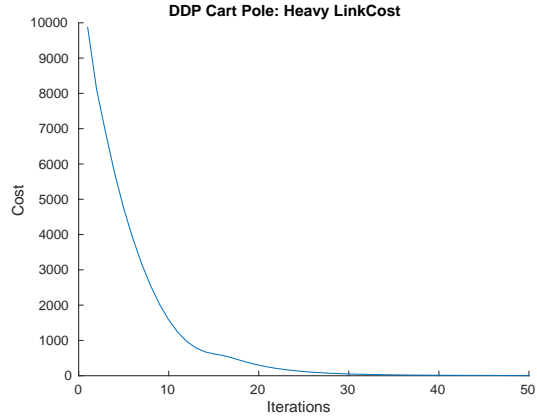}
	\caption{DDP Cart Pole Cost}
	\label{DCPCost}
\end{figure}		
 
We can view in Figure \ref{DCPControl} and \ref{GCPControl} the difference between the control effort for both algorithms. Interestingly enough, the pseudospectral algorithm came up with a lower final cost, even though it used higher maximum and minimum control. This observation is also related to the chosen masses for the Cart and Link. The link mass is much larger than the cart mass, hence creating a larger torque in the link by pushing the cart provides more power to swing up, rather than slowly pushing the cart side to side.

We can see in Figure \ref{DCPCost} the cost for DDP quickly converges by iteration 30. This will not always be the case, more complex systems, such as the Double Cart Pole will show spikes in the cost as DDP exhibits exploration into new control policies.

\subsection{Double Cart Pole}

The double cart pole is a highly under actuated system in which six states are controlled by a single control. It is similar to the cart pole system, however, there are now two links that must be controlled with just pushing the cart. Once again the task is to swing up both links from the stable equilibrium to the unstable equilibrium. Differential Dynamic Programming results in a feedforward, as well as a feedback gain, which proves key to solving this system. In order words, the system can apply full state feedback in order to reach the optimal trajectory. However, the Gauss Pseudospectral method has no such feedback term in its optimal control policy. As it result it is entirely feed forward exasperating the difficulty of the under actuated problem. This problem difficulty appears to manifest itself in the system runtime, GPOC takes 206 seconds to reach an error tolerance below 1E-6. DDP on the other hand can solve the problem in 59 seconds.

\begin{figure}[H]
	\centering
	\begin{subfigure}[b]{0.48\textwidth}
		\includegraphics[width=\textwidth]{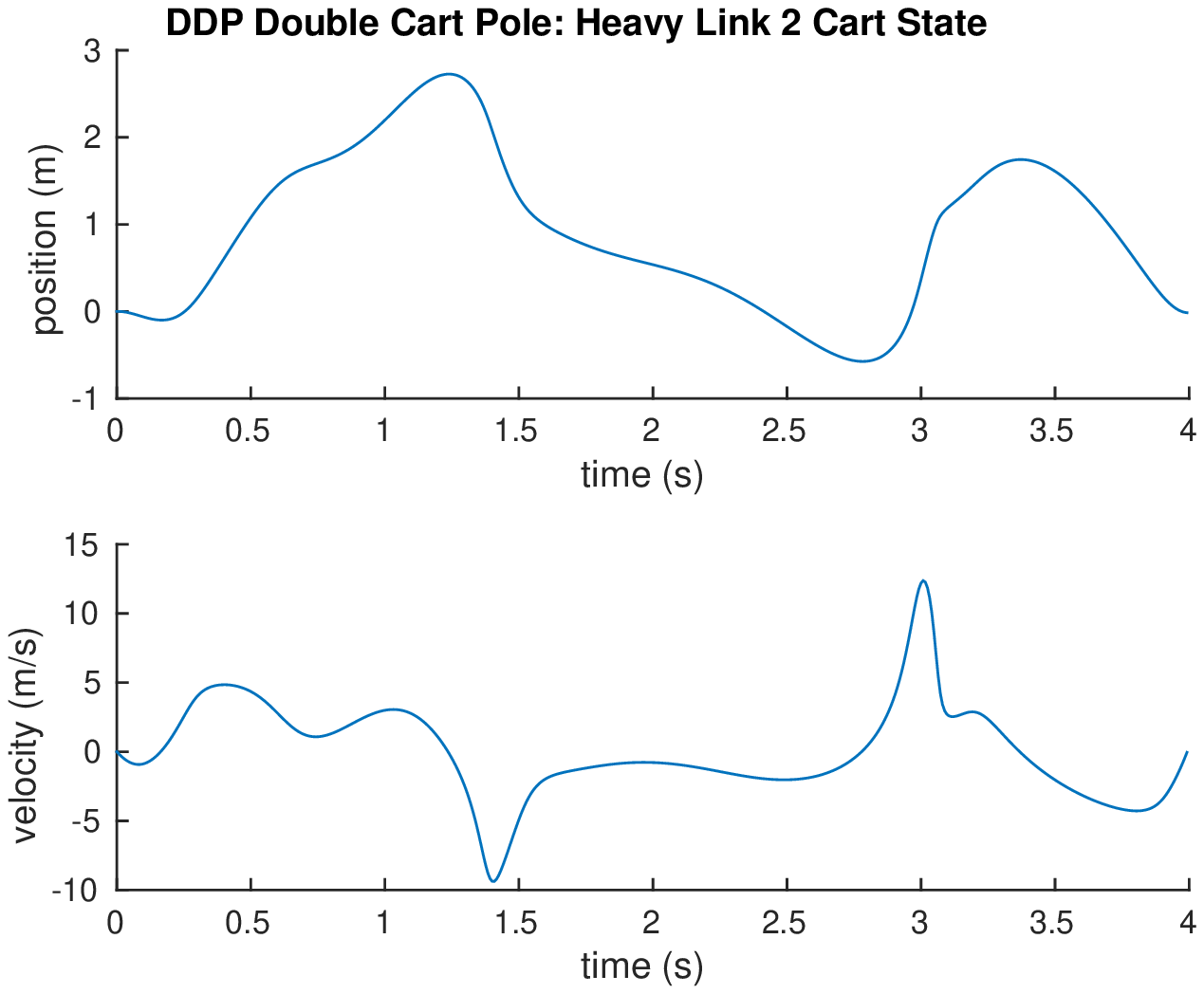}
		\caption{}
		\label{DDCPCS}
	\end{subfigure}
	~
	\begin{subfigure}[b]{0.48\textwidth}
		\includegraphics[width=\textwidth]{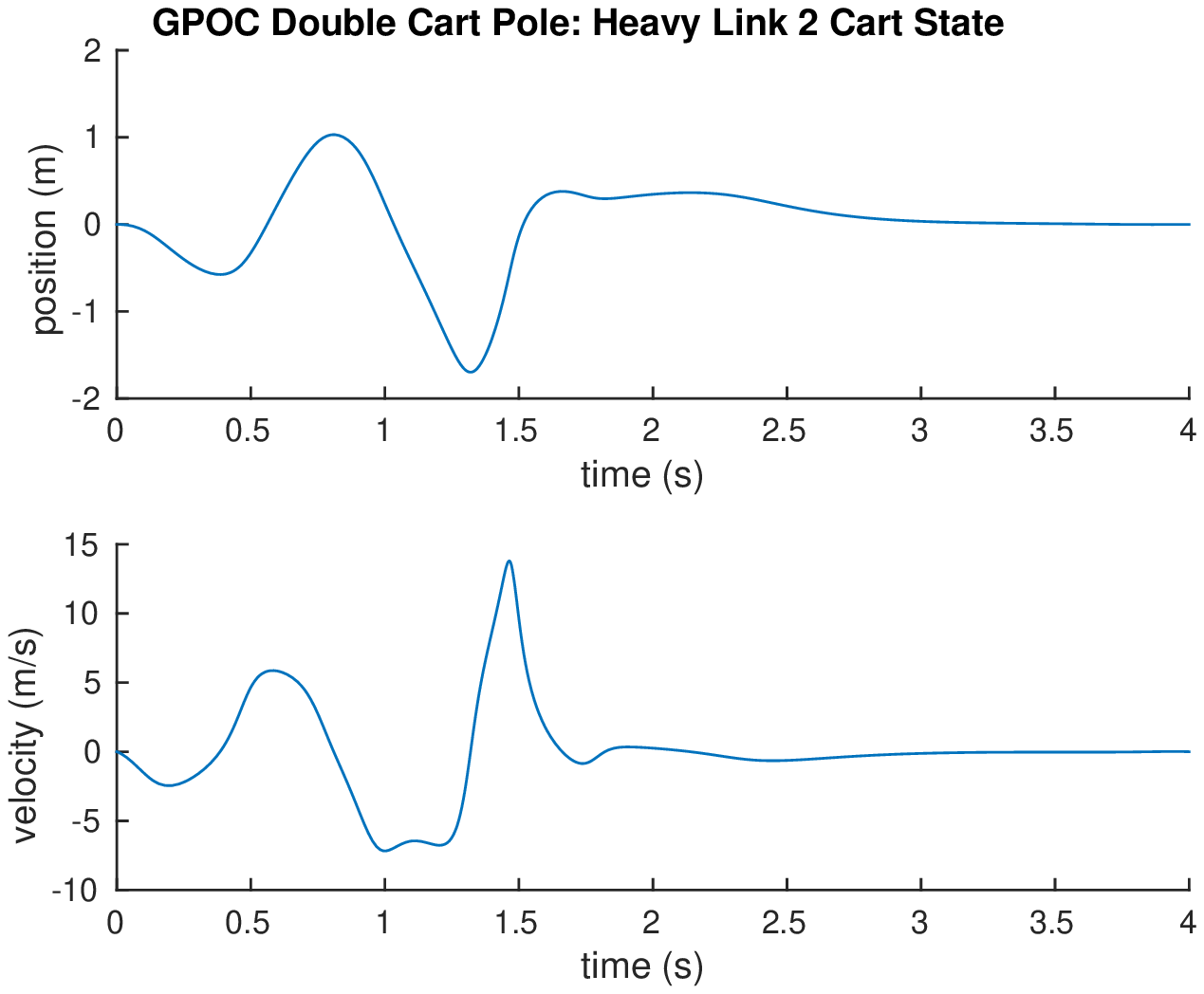}
		\caption{}
		\label{GDCPCS}
	\end{subfigure}
	\\
		\begin{subfigure}[b]{0.7\textwidth}
		\includegraphics[width=\textwidth]{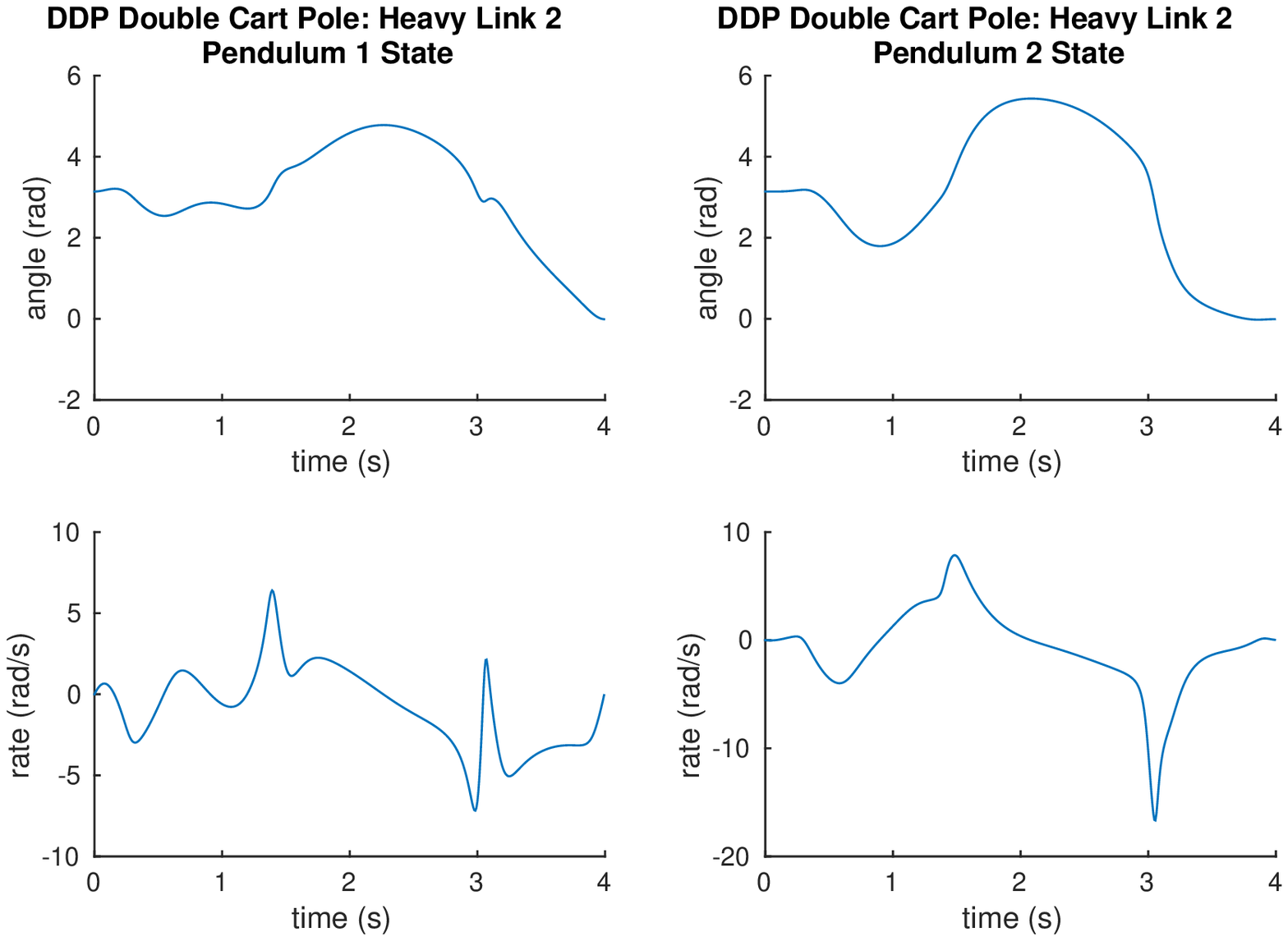}
		\caption{}
		\label{DDCPPS}
	\end{subfigure}
	\\
	\begin{subfigure}[b]{0.7\textwidth}
		\includegraphics[width=\textwidth]{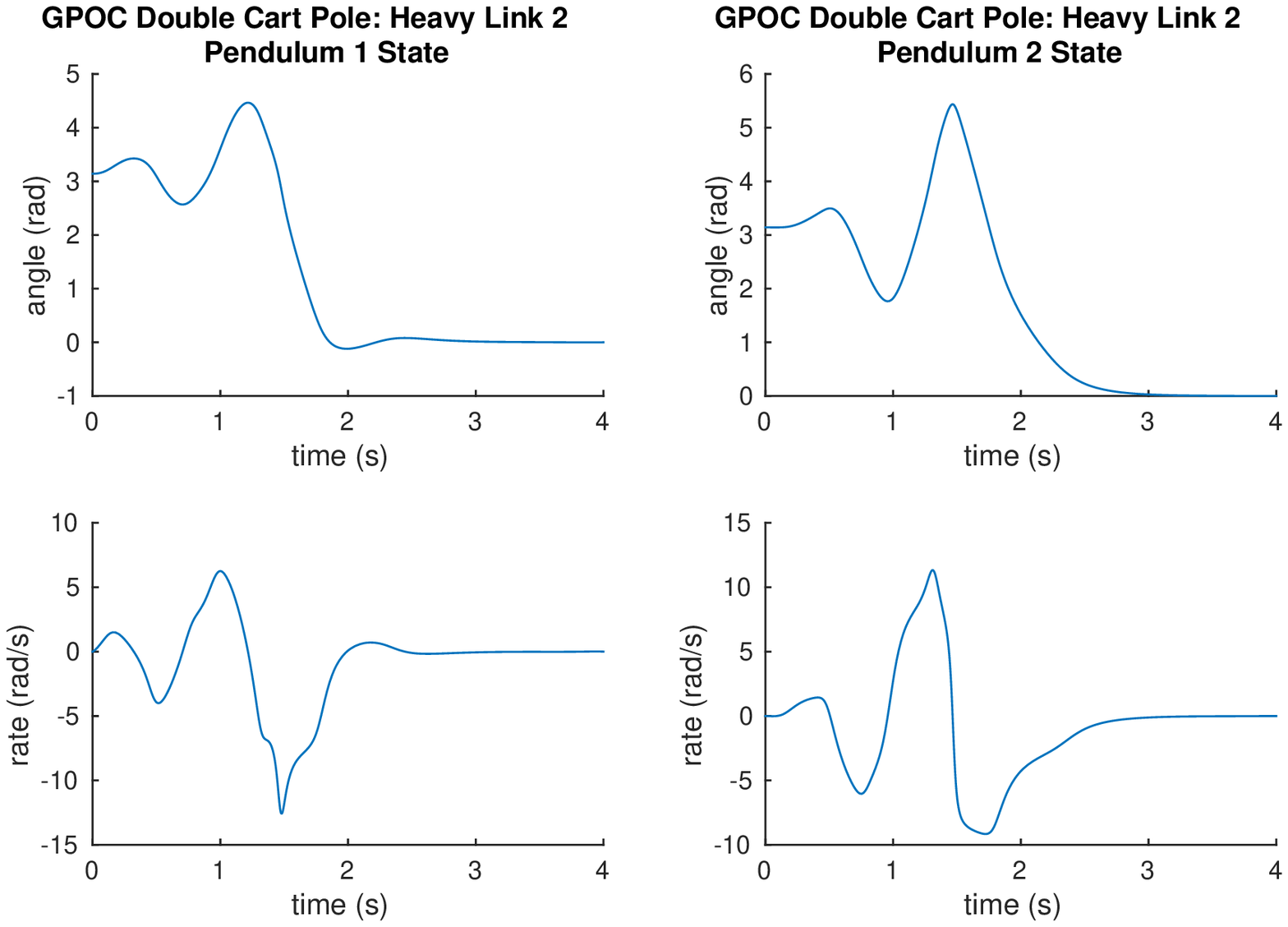}
		\caption{}
		\label{GDCPPS}
	\end{subfigure}
	\caption{Double Cart Pole State Comparison}
	\label{DCPoleCS}
\end{figure}

We can observe the optimal solution in Figure \ref{DCPoleCS}. The is a distinction in the trajectories that both methods provide. DDP take the entirety of the time horizon, all 4 seconds, to swing the two links up into the vertical position. We can see that in Figure \ref{GDCPCS}, the system manages to swing up in less time, and then stays near the origin for the remaining time. Once again, GPOC stays keeps the cart closer to the origin at the expense of higher maximum and minimum control effort.

\begin{figure}[H]
	\centering
	\begin{subfigure}[b]{0.48\textwidth}
		\includegraphics[width=\textwidth]{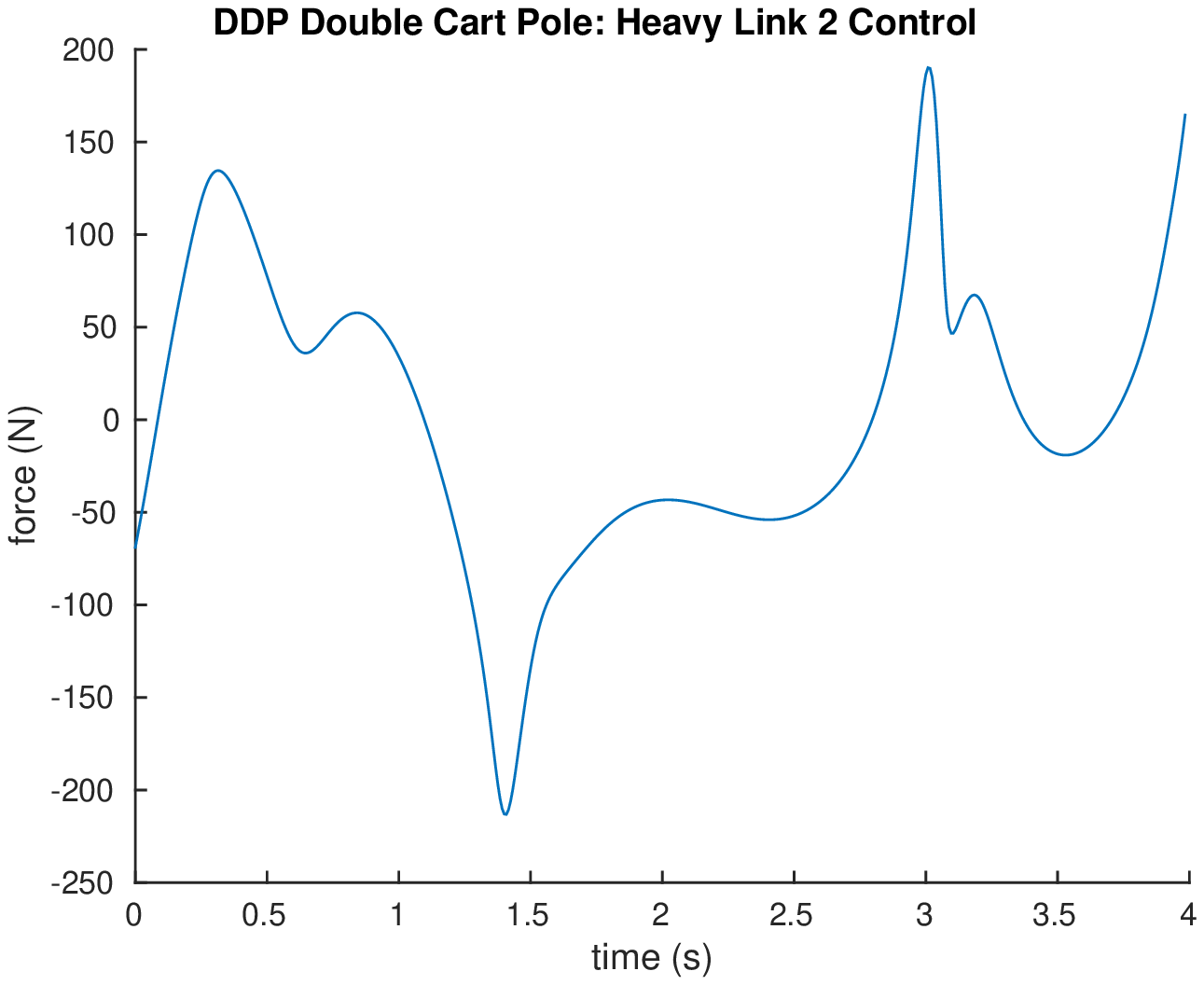}
		\caption{}
		\label{DDCPControl}
	\end{subfigure}
	~
	\begin{subfigure}[b]{0.48\textwidth}
		\includegraphics[width=\textwidth]{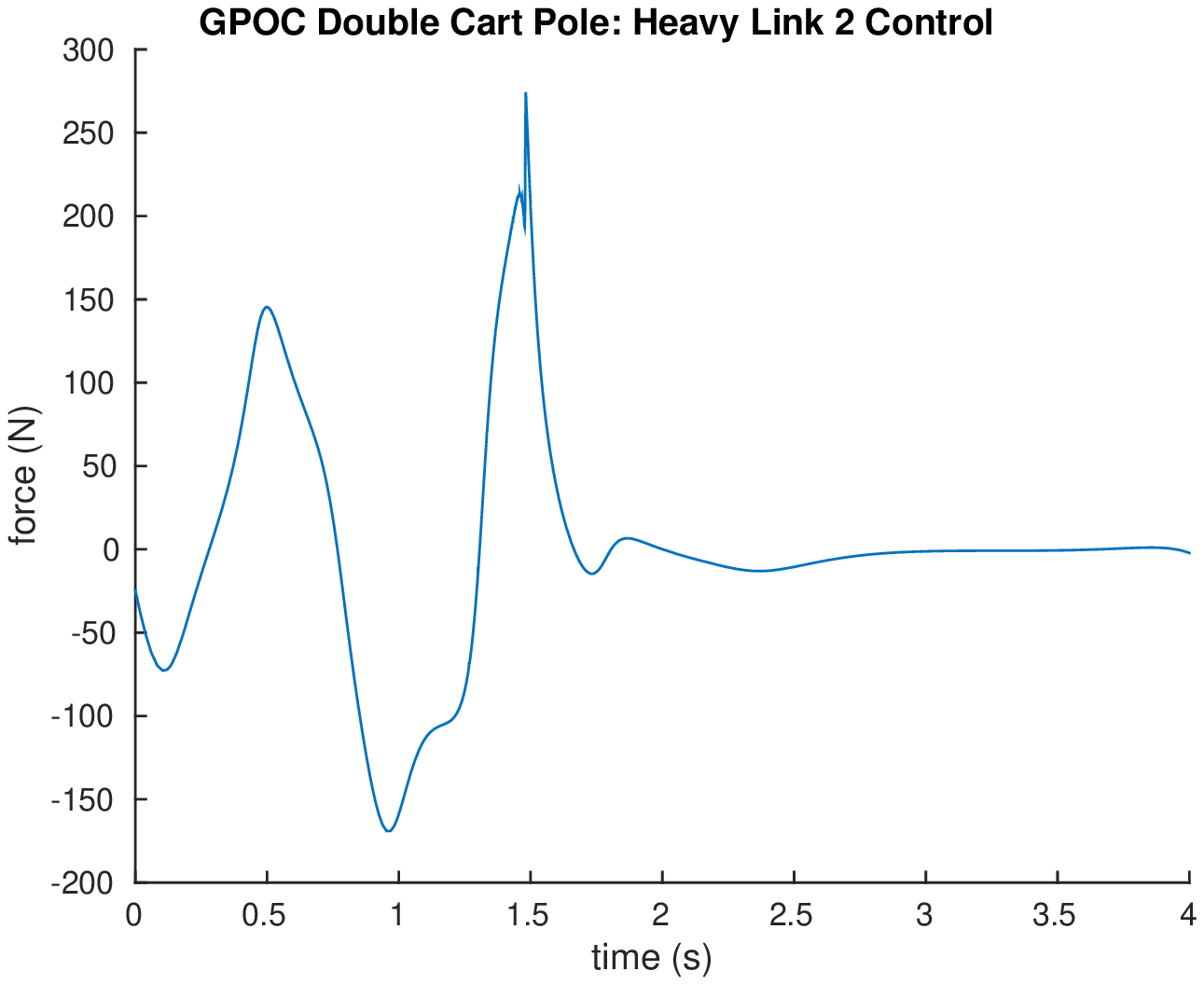}
		\caption{}
		\label{GDCPControl}
	\end{subfigure}
	\caption{Double Cart Pole Control Comparison}
	\label{DCPoleControl}
\end{figure}
	
\begin{figure}[H]
\centering
	\includegraphics[width=.48\textwidth]{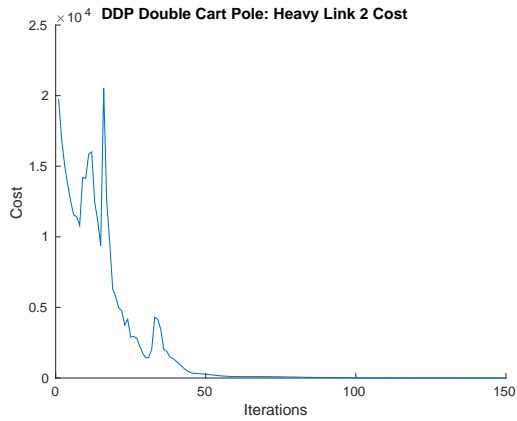}
	\caption{DDP Double Cart Pole Cost}
	\label{DDCPCost}
\end{figure}

Analyzing Figure \ref{DCPoleControl}, we can see the control effort differential manifest itself. Further enforcing GPOC's generic behavior, the control levels off to zero after 2.7 seconds, which is when the cart is almost to the origin and the links are almost vertical. The final cost of the DDP algorithm was .3632, while the final cost of the GPOC algorithm was 0.187. Even though the optimization took longer, the final solution was far more optimal, as long as control remains cheap.

We can see the oscillations in the cost for DDP in Figure \ref{DDCPCost}. As mentioned previously, DDP has the ability to explore for new policies given a particular quadratic structure of the cost function. Ultimately, the cost converges to a much more optimal minimum.

\subsection{Quadrotor}

Finally we come to the quadrotor system. The dynamical system has 12 states corresponding to position velocity, orientation, and angular rates. There are now four controls in this system as opposed to singular control in the other systems, these four controls correspond to the thrust of each rotor. The initial position was $[-1,1,.5]$ and the target state was $[.5,-1,1.5]$. The cost function now has an additional term which penalizes excessive orientation angles, such as the roll angle being 90 degrees.

\begin{figure}[H]
	\centering
	\begin{subfigure}[b]{0.41\textwidth}
		\includegraphics[width=\textwidth]{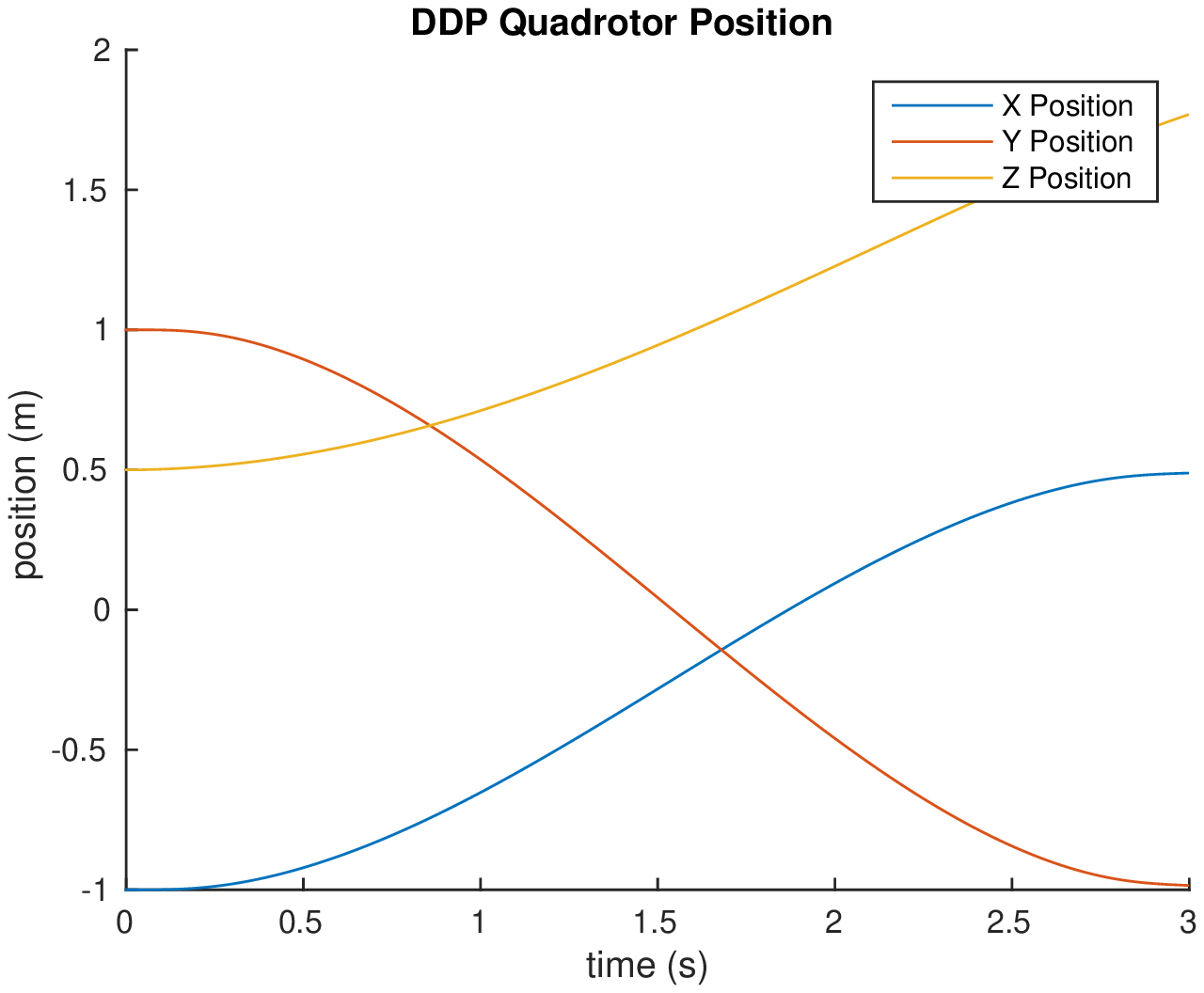}
		\caption{}
		\label{DQP}
	\end{subfigure}
	~
	\begin{subfigure}[b]{0.41\textwidth}
		\includegraphics[width=\textwidth]{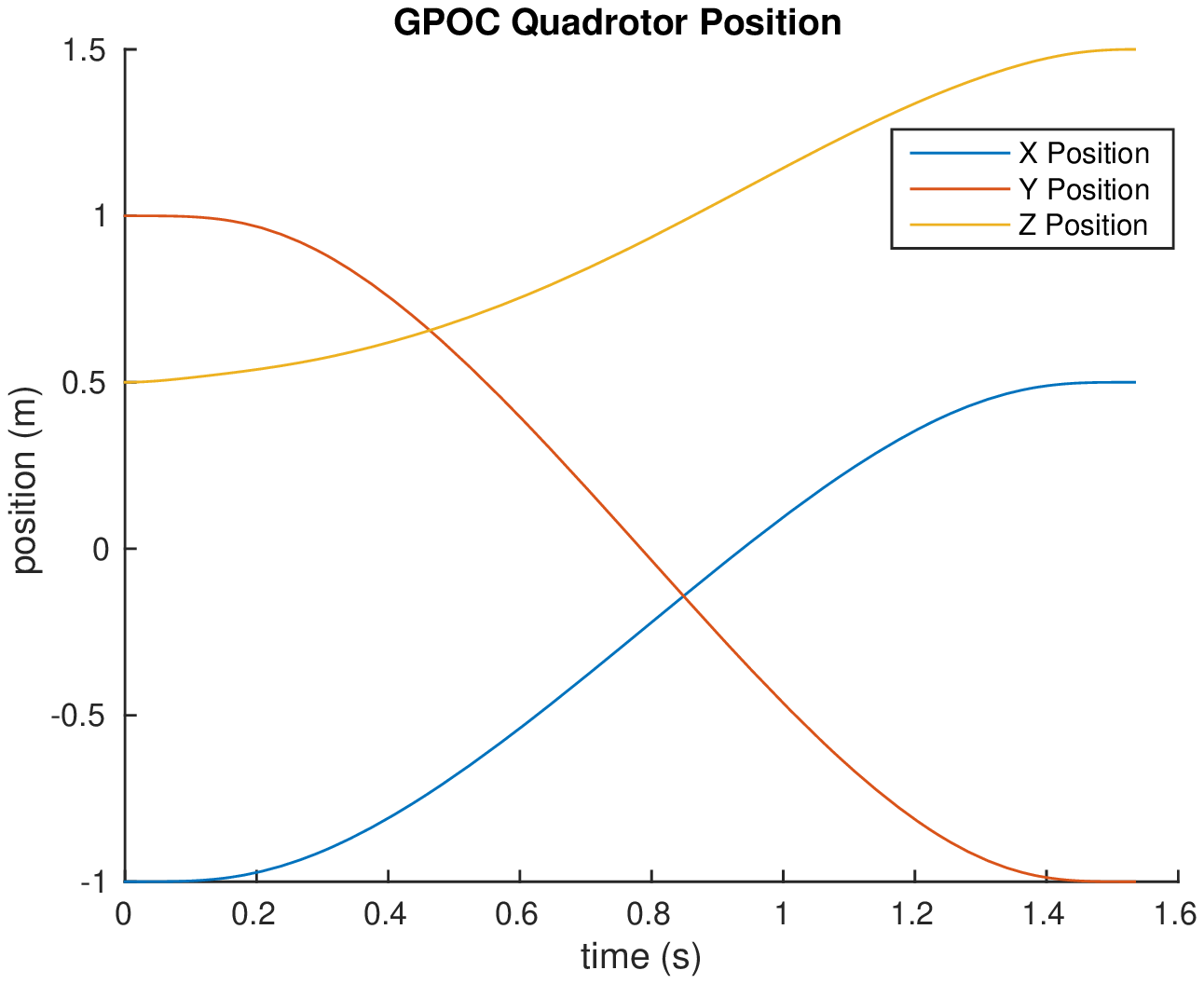}
		\caption{}
		\label{GQP}
	\end{subfigure}
	\\
	\begin{subfigure}[b]{0.41\textwidth}
		\includegraphics[width=\textwidth]{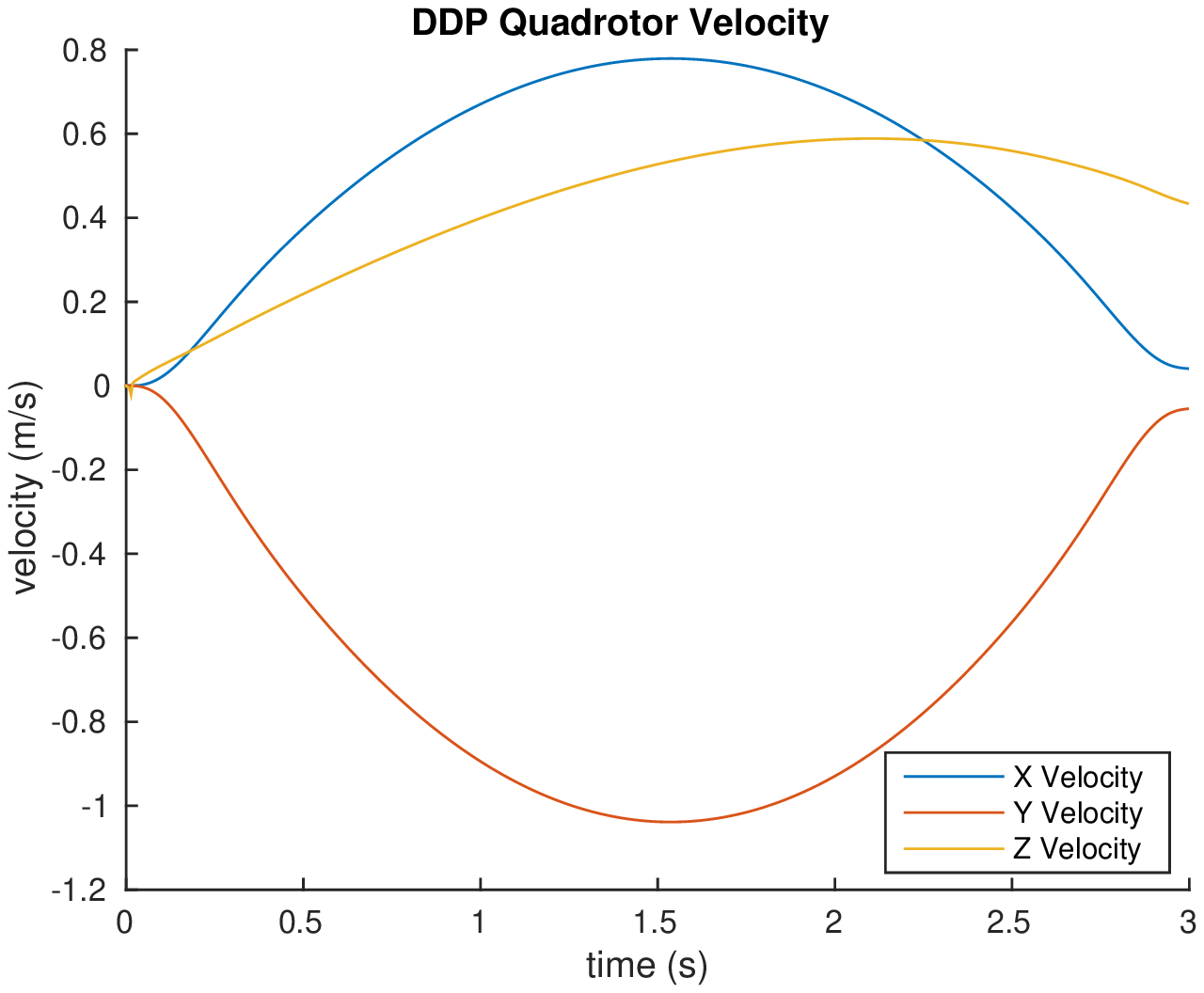}
		\caption{}
		\label{DQV}
	\end{subfigure}
	~
	\begin{subfigure}[b]{0.41\textwidth}
		\includegraphics[width=\textwidth]{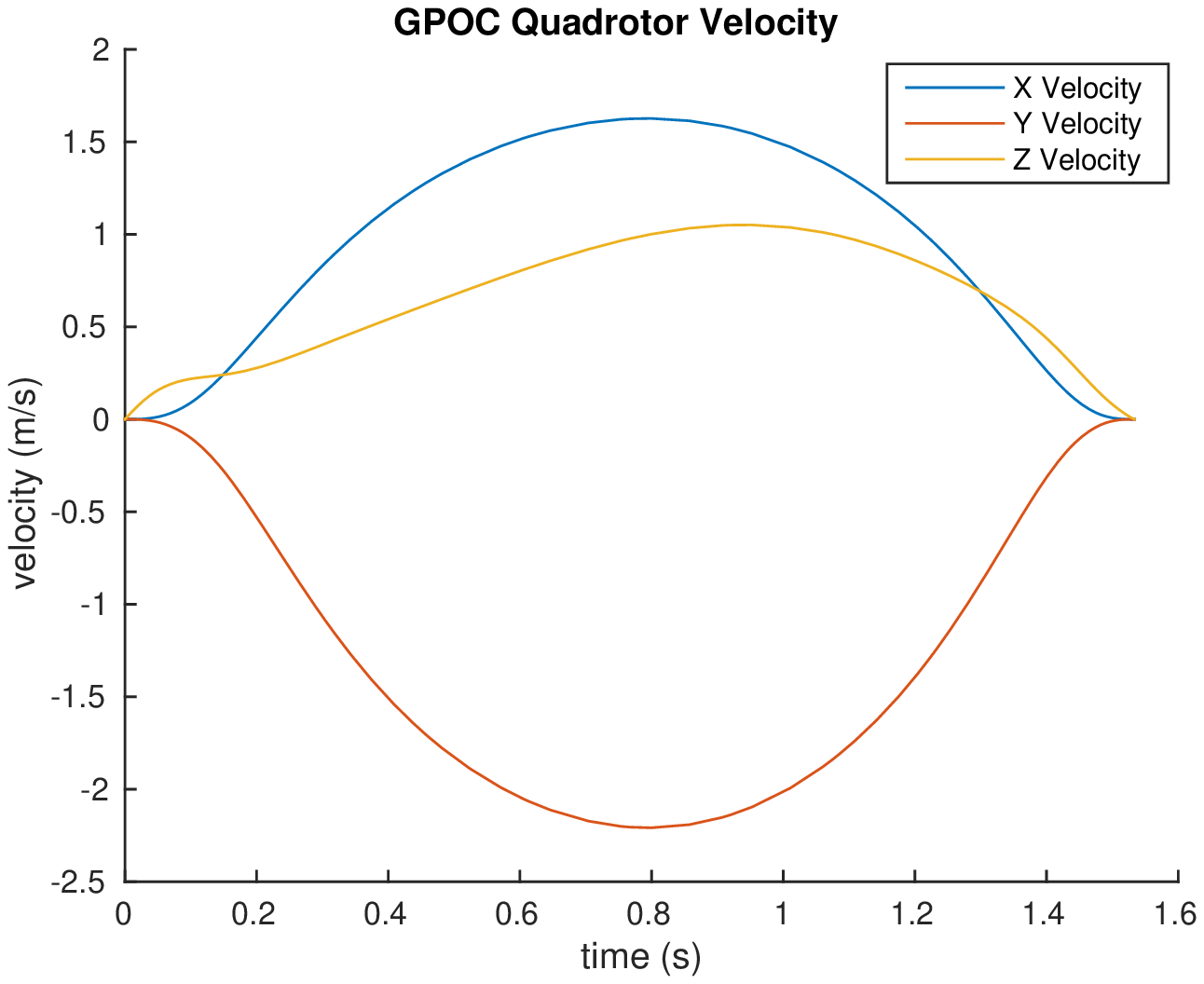}
		\caption{}
		\label{GQV}
	\end{subfigure}
	\begin{subfigure}[b]{0.41\textwidth}
		\includegraphics[width=\textwidth]{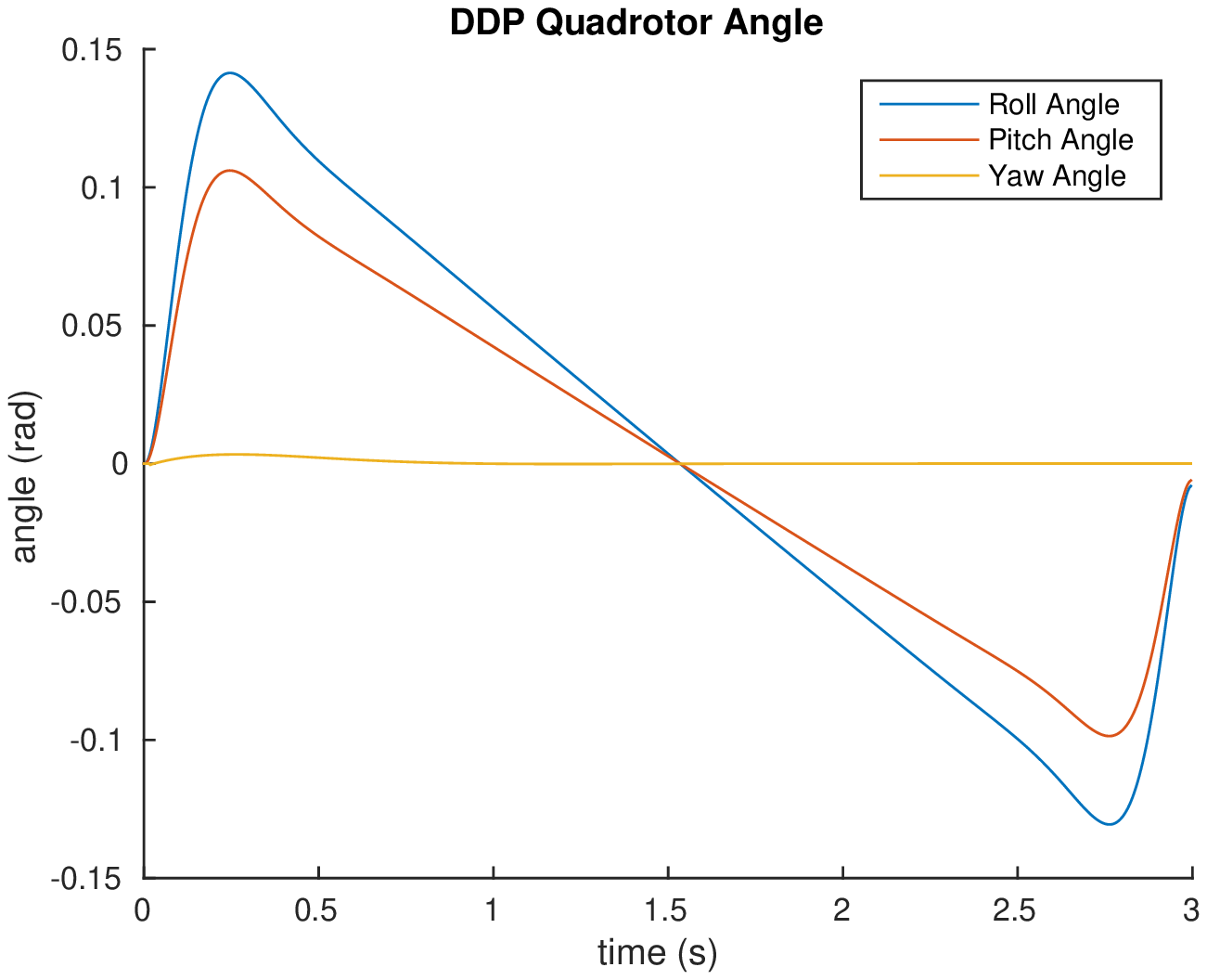}
		\caption{}
	\end{subfigure}
	~
	\begin{subfigure}[b]{0.41\textwidth}
		\includegraphics[width=\textwidth]{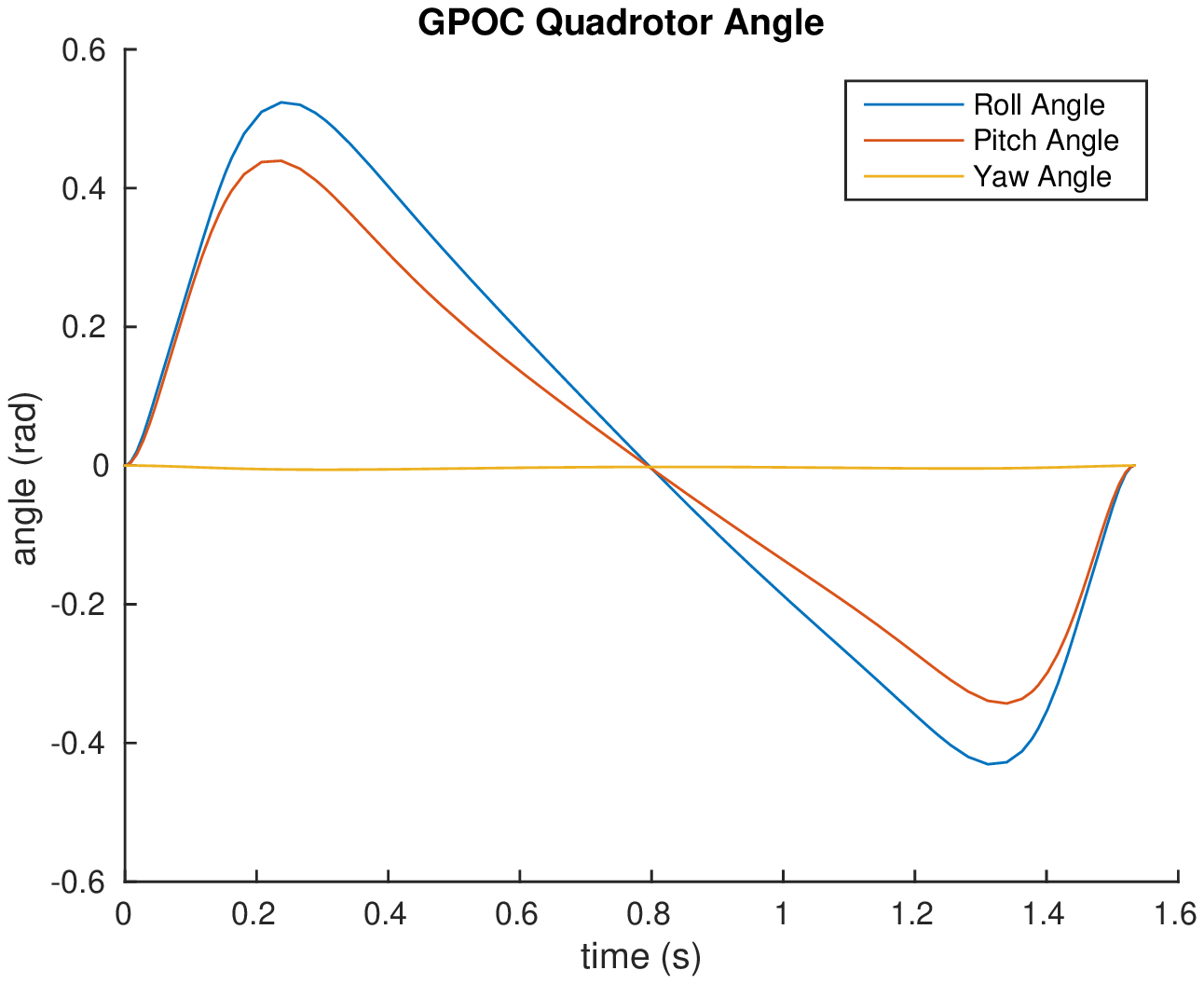}
		\caption{}
	\end{subfigure}
	\\
	\begin{subfigure}[b]{0.41\textwidth}
		\includegraphics[width=\textwidth]{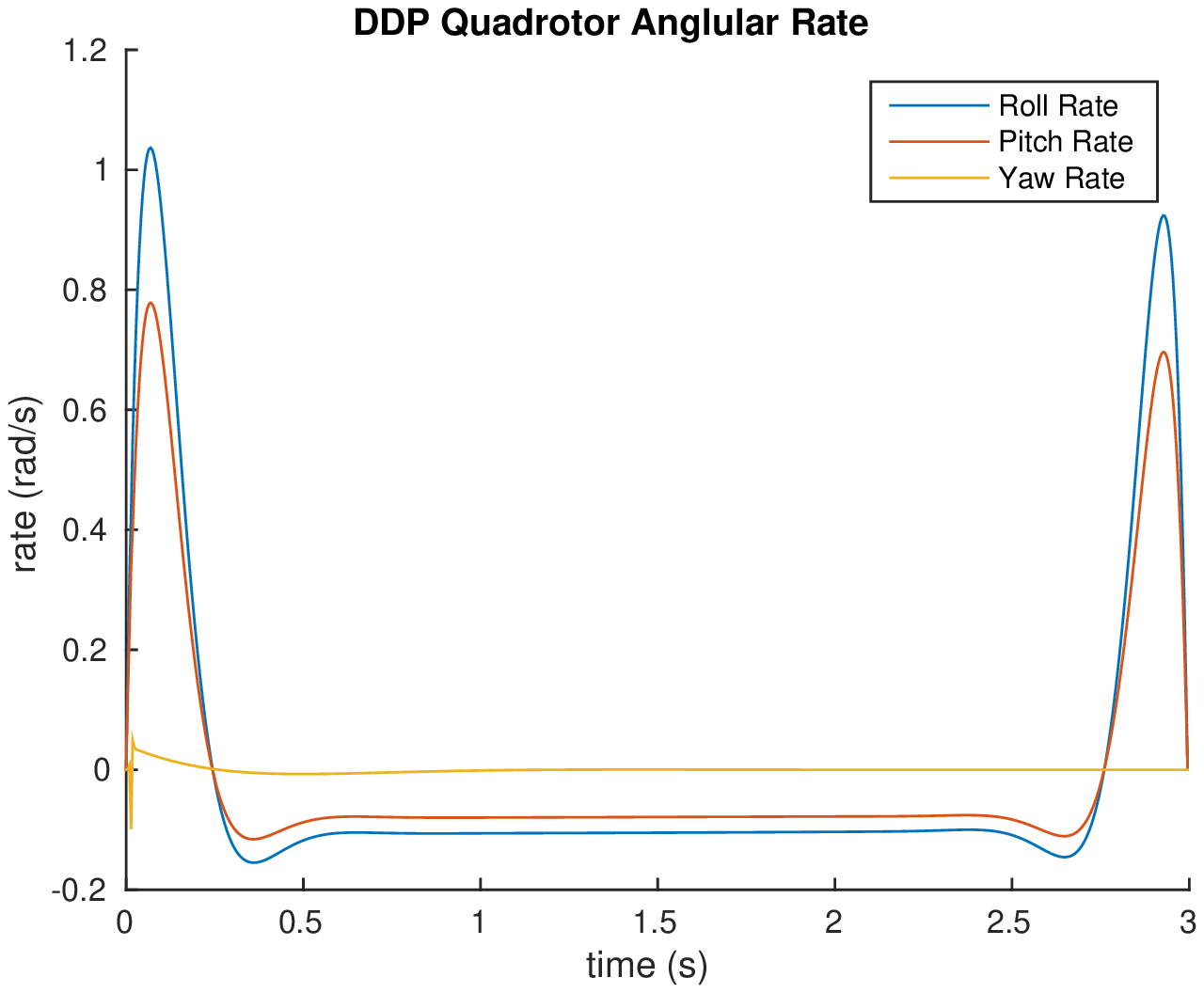}
		\caption{}
	\end{subfigure}
	~
	\begin{subfigure}[b]{0.45\textwidth}
		\includegraphics[width=\textwidth]{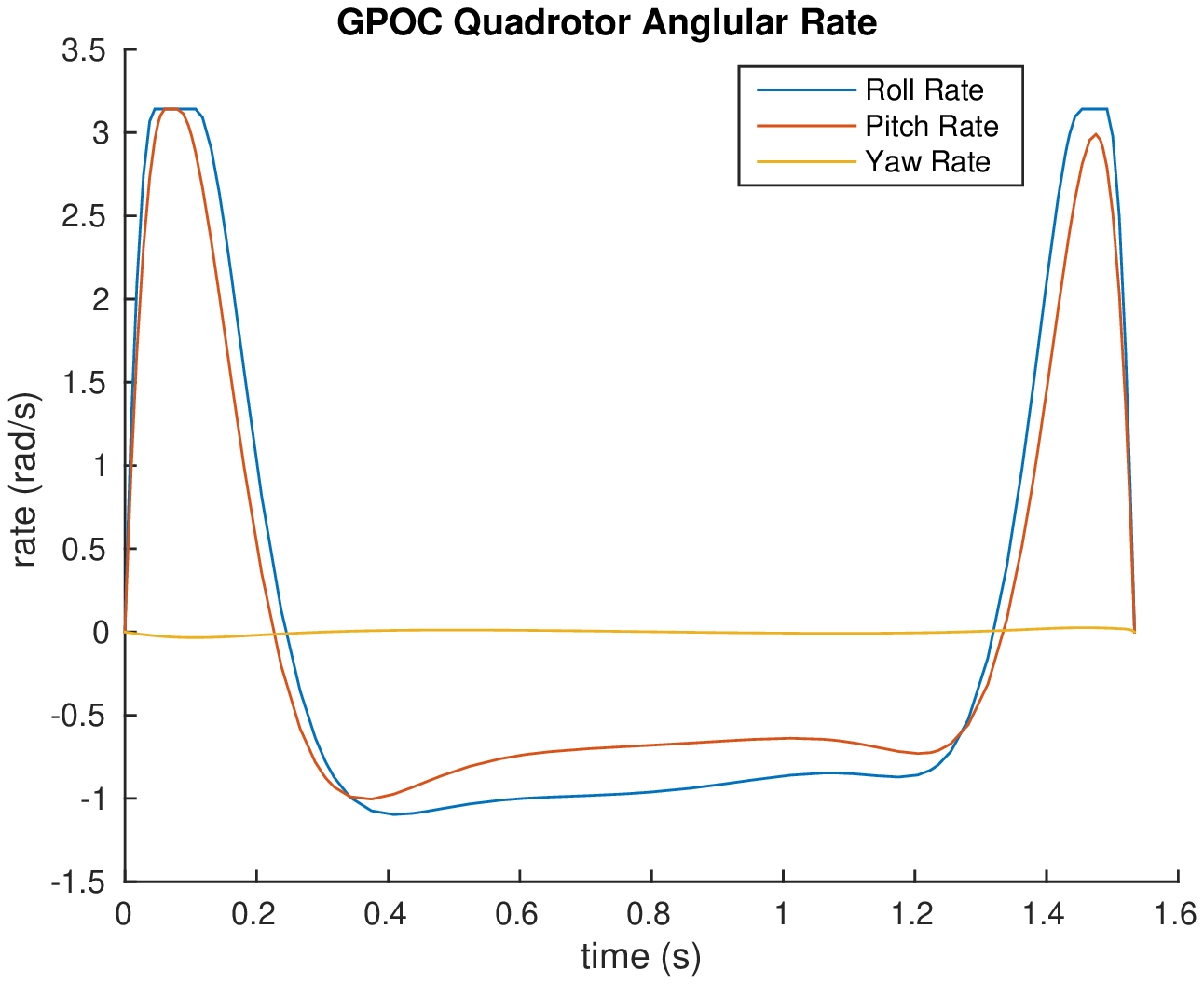}
		\caption{}
	\end{subfigure}
	\caption{Quadrotor State Comparison}
	\label{QuadState}
\end{figure}

\begin{figure}[H]
	\centering
	\begin{subfigure}[b]{0.48\textwidth}
		\includegraphics[width=\textwidth]{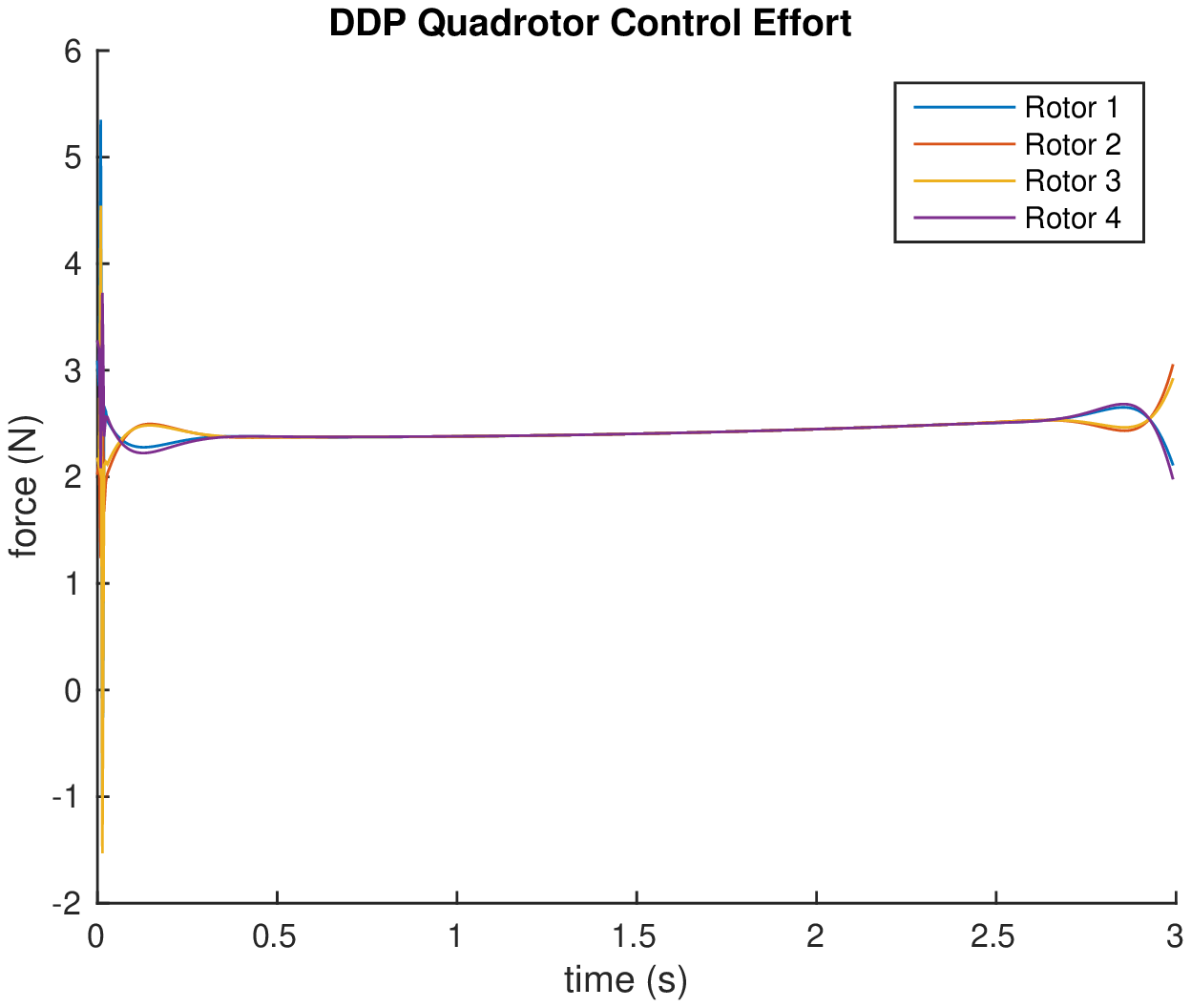}
		\caption{}
	\end{subfigure}
	~
	\begin{subfigure}[b]{0.48\textwidth}
		\includegraphics[width=\textwidth]{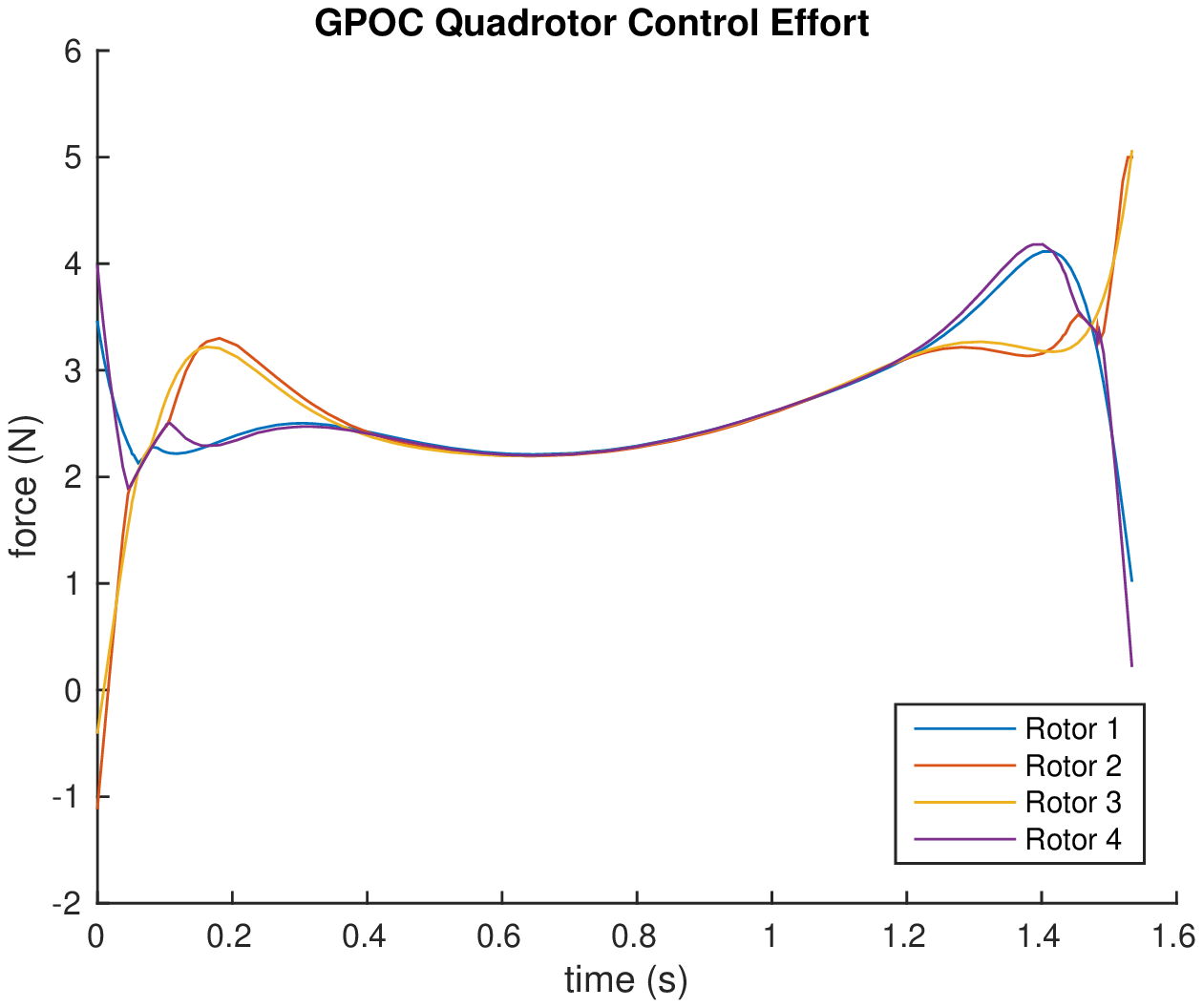}
		\caption{}
	\end{subfigure}
	\caption{Quadrotor Control Comparison}
	\label{QuadControl}
\end{figure}
	
\begin{figure}[H]
\centering
	\includegraphics[width=.48\textwidth]{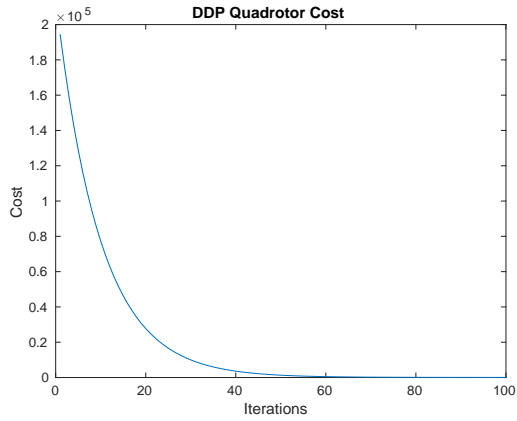}
	\caption{DDP Quadrotor Pole Cost}
\end{figure}

We can observe the same patterns between the two algorithms. GPOC utilizes more force, but reaches the target faster, while DDP takes the entire time horizon. Note that DDP was extremely sensitive to the discrete time step size, as well as the choice of weights in the cost function. DDP converged with a cost of 44.6182 with a runtime of 41 seconds. GPOC converged with a cost of 47.9891 and with a runtime of 18 seconds. Once again GPOC is faster, and notice in Figure \ref{QuadState}, DDP did not converge precisely to the desired target state, thus GPOC delivered a more optimal solution. This time we see that angle, velocity, and angular rate were higher with GPOC as opposed to DDP. This is likely related to the increased amount of control seen in Figure \ref{QuadControl}, but also to how sensitive GPOC and DDP are to weights in the cost function. For GPOC, the final state must exactly match the bounds described in the problem, while for DDP, the optimal trajectory need not necessary reach the final target exactly, again depending on how heavily the weights are set. This is a fundamental difference in the behavior of each algorithm. Finally, we can observe that this rotor configuration is standard and does not reflect on how well each algorithm withstands crazy configurations, such as large rotor to CoM lengths, or non-symmetric inertias. This was explored partially by the previous experiments, as the carts were made to be far lighter than the links, to increase the difficulty of optimization.

\section{Future Work}

These simulations are just the beginning for the analysis of these two trajectory optimization techniques. Ultimately the goal is to apply both methods into an ROS simulation of the Barret Technologies WAM arm and our quadrotor. In order to do this, the Discrete Time DDP Implementation must be completed in C++, and it would also be best to analyze the performance of the continuous time formulation against GPOC. In addition, we want to extend the capabilities of these algorithms, namely into dynamics with contact and stochasticity. Currently, I am studying methods to propagate uncertainty through polynomial chaos, as well as methods to solve a linear complementary problem with DDP to incorporate discontinuity. Finally, both methods can also be studied with a receding horizon framework, to assess their ability to be applied to real systems. DDP has the advantage of feedback gains pushing the system towards the optimal trajectory, GPOC has the ability to divide the trajectory into small bits, to avoid the need for feedback by optimizing more often and faster. There is much to be done in the future, but the possibilities are promising.

\newpage
\bibliographystyle{plain}
\bibliography{trbib}
\end{document}